\pdfoutput=1
\RequirePackage{ifpdf}
\ifpdf 
\documentclass[pdftex]{sigma}
\else
\documentclass{sigma}
\fi

\newcommand{\R}{\mathbb{R}}
\newcommand{\h}{\mathbb{H}}
\newcommand{\ei}{\varepsilon_i}
\newcommand{\ej}{\varepsilon_j}
\newcommand{\ek}{\varepsilon_k}

\numberwithin{equation}{section}

\begin{document}

\allowdisplaybreaks

\renewcommand{\thefootnote}{$\star$}

\renewcommand{\PaperNumber}{033}

\FirstPageHeading

\ShortArticleName{The Symmetry Group of Lam\'e's System and  Conformally Flat Hypersurfaces}

\ArticleName{The Symmetry Group of Lam\'e's System\\ and the Associated Guichard Nets\\ for Conformally Flat Hypersurfaces\footnote{This paper is a~contribution to the Special Issue  ``Symmetries of Dif\/ferential
Equations:  Frames, Invariants and~Applications''.
The full collection is available at
\href{http://www.emis.de/journals/SIGMA/SDE2012.html}{http://www.emis.de/journals/SIGMA/SDE2012.html}}}

\Author{Jo\~{a}o Paulo dos SANTOS and Keti TENENBLAT}
\AuthorNameForHeading{J.P.~dos Santos and K.~Tenenblat}
\Address{Departamento de Matem\'{a}tica, Universidade de Bras\'{i}lia, 70910-900, Bras\'{i}lia-DF, Brazil}
\Email{\href{mailto:j.p.santos@mat.unb.br}{j.p.santos@mat.unb.br}, \href{mailto:k.tenenblat@mat.unb.br}{k.tenenblat@mat.unb.br}}

\ArticleDates{Received October 01, 2012, in f\/inal form April 12, 2013; Published online April 17, 2013}

\Abstract{We consider conformally f\/lat hypersurfaces in four dimensional space forms with their associated Guichard nets and Lam\'e's system of equations. We show that the symmetry group of the Lam\'e's system, satisfying  Guichard condition, is given by translations and dilations in the independent variables and dilations in the dependents variables. We obtain the solutions which are invariant under the action of the 2-dimensional subgroups of the symmetry group. For the solutions which are invariant under translations, we obtain the corresponding conformally f\/lat hypersurfaces and we describe the corresponding Guichard nets. We show that the coordinate surfaces of the Guichard nets have constant Gaussian curvature, and the sum of the three curvatures is equal to zero. Moreover, the Guichard nets are foliated by f\/lat surfaces with constant mean curvature. We prove that there are solutions of the Lam\'e's system, given in terms of Jacobi elliptic functions, which are invariant under translations, that correspond to a new class of conformally f\/lat hypersurfaces.}


\Classification{53A35; 53C42}

\renewcommand{\thefootnote}{\arabic{footnote}}
\setcounter{footnote}{0}

\section{Introduction}

The investigation of conformally f\/lat hypersurfaces has been of interest for quite some time. Any surface in $\R^3$ is conformally  f\/lat, since it can be parametrized by isothermal coordinates. For higher dimensional hypersurfaces, E. Cartan \cite{Cartan} gave a complete classif\/ication for the conformally f\/lat hypersurfaces of an $(n+1)$-dimensional space form when $n+1 \geq 5$. He proved that such hypersurfaces are quasi-umbilic, i.e., one of the principal curvatures has multiplicity at least $n-1$. In the same paper, Cartan investigated the case $n + 1 = 4$ . He showed that the quasi-umbilic surfaces are conformally f\/lat, but the converse does not hold (for a proof see \cite{Lafontaine}). Moreover, he gave a characterization of the conformally f\/lat 3-dimensional hypersurfaces, with three distinct principal curvatures, in terms of certain integrable distributions. Since then, there has been an ef\/fort to obtain a classif\/ication of hypersurfaces satisfying Cartan's characterization.

Lafontaine \cite{Lafontaine} considered hypersurfaces of type $M^3 = M^2 \times I \subset \R^4$. He obtained the following classes of conformally f\/lat hypersurfaces: a)~$M^3$ is a cylinder over a surface,  $M^2 \subset \R^3$, with constant curvature;
b)~$M^3$ is a cone over a surface in the sphere,  $M^2 \subset \mathbb{S}^3$, with constant curvature;
c)~$M^3$ is obtained by rotating a constant curvature surface of the hyperbolic space, $M^2 \subset \h^3 \subset \R^4$, where $\h^3$ is the half space model.

Motivated by Cartan's paper, Hertrich-Jeromin \cite{Jeromin1}, established a correspondence between conformally f\/lat three-dimensional hypersurfaces, with three distinct principal curvatures,  and Guichard nets. These are systems of triply orthogonal surfaces originally considered by C.~Gui\-chard in \cite{Guichard}, where he referred  to those systems as the analogues of isothermal coordinates.

In view of Hertrich-Jeromin results, the problem of classifying conformally f\/lat 3-dimensional hypersurfaces was transferred to the problem of classifying Guichard nets in $\R^3$. These are open sets of $\R^3$, with an orthogonal f\/lat metric $g= \sum\limits_{i=1}^3 l_i^2 dx_i^2$, where the functions $l_i$ satisfy the Guichard condition, namely,
\begin{gather*}
l_1^2 - l_2^2 + l_3^2 = 0,
\end{gather*}
and a system of second-order partial dif\/ferential equations, which is called Lam\'e's system (see~(\ref{lame})).

Hertrich-Jeromin obtained an example of a Guichard net, starting from surfaces parallel to Dini's helix and he proved that the corresponding conformally f\/lat hypersurface was a new example, since it did not belong to the class described by Lafontaine.

In \cite{Suyama1,Suyama2, Suyama3}, Suyama extended the previous results by showing that the Guichard nets described by Hertrich-Jeromin are characterized in terms of a dif\/ferentiable function $\varphi(x_1, x_2, x_3)$ that determines, up to conformal  equivalence, the f\/irst and second fundamental forms of the corresponding conformally f\/lat hypersurfaces. Moreover, Suyama showed that if $\varphi$ does not depend on one of the variables, then the hypersurface is conformal to one of the classes described by Lafontaine. He also showed that the function associated to the  example given by Hertrich-Jeromin satisf\/ied $\varphi_{,x_1 x_2} = \varphi_{,x_2 x_3} = 0.$ Starting with this condition on $\varphi$, Suyama obtained a partial classif\/ication of such conformally f\/lat hypersurfaces. The complete classif\/ication of conformally f\/lat hypersurfaces, satisfying the above condition on the partial derivatives of $\varphi$, was obtained by Hertrich-Jeromin and Suyama in \cite{suyamajer}. They showed that these hypersurfaces correspond to a special type of Guichard nets. The authors called them cyclic Guichard nets, due to the fact that one of the coordinates curves is contained in a circle.

In this paper, we obtain solutions $l_i$ satisfying Lam\'e's system and the Guichard condition, which are invariant under the action of the 2-dimensional subgroups of the symmetry group of the system. Moreover, we investigate the properties of the Guichard nets and of the conformally f\/lat hypersurfaces associated to the solutions $l_i$. We f\/irst determine the symmetry group of Lam\'e's system satisfying the Guichard condition. We prove that the group is given by translations and dilations of the independent variables $x_i$ and dilations of the dependent variables $l_i$.

We obtain the solutions $l_i$, $i =1, 2, 3$, which are invariant under the action of the 2-dimensional translation subgroup, i.e., $l_i (\xi)$, where $\xi = \sum\limits_{i=1}^3 \alpha_i x_i.$ These solutions are given explicitly in Theorem~\ref{teo.nonconstant.trans} by Jacobi elliptic functions, whenever all the functions $l_i$ are not constant and in Theorem~\ref{uma.constante} when one of the functions $l_i$ is constant. Moreover, we consider the solutions $l_i$ which are invariant under the 2-dimensional subgroup involving translations and dilations, i.e., $l_i (\eta)$, where $\eta = \sum\limits_{j=1}^3 a_j x_j/ \sum\limits_{k=1}^3 b_k x_k$. In this case, if we require the functions $l_i(\eta)$ to depend on all three variables, then $l_i$ are constant functions. Otherwise, the solutions~$l_i(\eta)$ are given explicitly in Theorem~\ref{uma.constante.dilata}. The symmetry subgroup of dilations on the dependent variables is irrelevant for the study of conformally f\/lat hypersurfaces.

Considering the functions $l_i$ which are invariant under the action of translations, we study the corresponding Guichard nets. We show that their coordinate surfaces have constant Gaussian curvature and the sum of the three curvatures is equal to zero. Moreover the Guichard nets are foliated by f\/lat surfaces, with constant mean curvature.

Finally, we investigate the conformally f\/lat hypersurfaces associated to the functions $l_i$ which are invariant under the action of translations. We show that, whenever the basic invariant~$\xi$ depends on two variables, the hypersurface is conformal to one of the products considered by Lafontaine. In this case, the three-dimensional conformally f\/lat hypersurfaces are constructed from f\/lat surfaces contained in the hyperbolic 3-space~$\h^3$ or in the sphere $\mathbb{S}^3$. Whenever the basic invariant $\xi$ depends on all three independent variables, then the functions~$l_i (\xi)$, which are given in terms of Jacobi elliptic functions, produce a new class of conformally f\/lat hypersurfaces.

In Section~\ref{section2}, we review the correspondence between conformally f\/lat 3-dimensional hypersurfaces with Lam\'e's system, and Guichard nets.

In Section~\ref{section3}, we obtain the symmetry group of Lam\'e's system satisfying Guichard condition and the solutions which are invariant under 2-dimensional subgroups of the symmetry group. The motivation and the technique used in this section were inspired by the fact that our system of dif\/ferential equations is quite similar to the intrinsic generalized wave and sine-Gordon equations and the generalized Laplace and sinh-Gordon equations. The symmetry groups of these systems and the solutions invariant under subgroups were obtained by Tenenblat and Winternitz in~\cite{keti3} and Ferreira~\cite{ferreira}.
The geometric properties of the submanifolds corresponding to the solutions invariant under the subgroups of symmetries can be found in~\cite{barbosa} and~\cite{rabelo}.

In Sections~\ref{section4} and~\ref{section5}, we describe the geometric properties of the Guichard nets and of the conformally f\/lat hypersurfaces that are associated to the solutions of Lam\'e's system which are invariant under the action of the translation group.

The solutions $l_i$ of Lam\'e's system, satisfying Guichard condition, which are invariant under the subgroup of dilations of the independent variables and the corresponding geometric theory, will be considered in another paper. Such solutions are obtained by solving a (reduced) system of partial dif\/ferential equations, in contrast to what occurs in this paper, where  the Lam\'e's system is reduced to a system of ordinary dif\/ferential equations.

\section{Lam\'e's system and conformally f\/lat hypersurfaces}\label{section2}

Consider the Minkowski space $\R^6_1$ with coordinates $(x_0, \ldots, x_5)$ and the scalar product $\langle \; , \; \rangle $ given by
\begin{gather*}
\langle \; ,\; \rangle    : \   \R^6 \times \R^6   \longrightarrow   \R ,\\
\hphantom{\langle \; ,\; \rangle    : \ } \ (v,w)   \mapsto   -v_0w_0 +   \sum_{i=1}^5 v_iw_i.
\end{gather*}
Let $L^5\! =\! \left\{ y \in \R^6_1 \,|\, \langle  y, y \rangle\!  =\!0 \right\}$, be the light cone in~$\R^6_1$ and consider $m_K \in \R^6_1$, with $\langle  m_K, m_K \rangle  \!=\! K$. Then, it is not dif\/f\/icult to see that, the sets
\begin{gather*}
	M^4_K = \left\{ y \in L^5 \,|\, \langle  y, m_K \rangle  = -1 \right\}, 
\end{gather*}
with the metric induced from $\R^6_1$, are complete Riemannian manifolds with constant sectional curvature~$K$. If $K<0$, then $M^4_K$ consists of two connected components which can be isometrically identif\/ied (see \cite[Lemma~1.4.1]{Jeromin2} for details).

With this approach, consider a Riemannian immersion $f: M^3 \rightarrow M^4_K \subset L^5$, with unit normal~$n$. Then  $\langle  df, n \rangle  \equiv 0$, and $n$ also satisf\/ies $\langle  n, m_K \rangle  = \langle  n, f \rangle  = 0$. Let  $\tilde{f} : M^3 \rightarrow L^5$ be an immersion given by $\tilde{f} = e^u f$, where $u$ is a dif\/ferentiable function on $M$. Observe that the metric induced on $\tilde{f}$ is conformal to the metric induced on the immersion $f$, i.e.,
\begin{gather*}
\langle  d \tilde{f}, d \tilde{f} \rangle  = e^{2u} \langle  d f, d f \rangle .
\end{gather*}

\begin{definition}
Let $f:M^3 \rightarrow L^5$ be an immersion such that the induced metric, $\langle  d f, d f \rangle $, is positive def\/inite. Let $n$ be  a unit normal with $\langle  f,n \rangle =0$ and consider dif\/ferentiable functions $u$ and $a$ on $M^3$. Then the pair $(f,n)$ is called a \textit{strip} and the pair $(\tilde{f}, \tilde{n})$ given by{\samepage
\begin{gather*}
 \tilde{f} = e^uf  , \qquad \tilde{n} = n+ af
\end{gather*}
  is called a \textit{conformal deformation} of the strip $(f,n)$.}
\end{definition}

Therefore, we can deform a conformally f\/lat immersion in a space form $f:M^3 \rightarrow M^4_K \subset L^5$ to a f\/lat immersion in the light cone $\tilde{f}:M^3 \rightarrow L^5$, by considering a conformal deformation, and vice-versa.
Hence the problem of investigating conformally f\/lat hypersurfaces in space forms reduces to a problem of studying  f\/lat immersions in the light cone $f: M^3 \rightarrow L^5$. We say that a conformally f\/lat hypersurface in a space form $M^4_k$ is \emph{generic} if it has three distinct principal curvatures.  Hertrich-Jeromin in \cite{Jeromin1} established a relation between generic conformally f\/lat hypersurfaces in $M^4_k$ and Guichard nets \cite{Guichard}. Namely, let $e_1$, $e_2$, $e_3$ be an orthonormal frame tangent to $M^3 \subset M^4_k$, such that $e_i$ are principal directions. Let $\omega_1$, $\omega_2$, $\omega_3$ be the co-frame and let $k_1$, $k_2$, $k_3$ be the principal curvatures. Assume that locally $k_3 > k_2 > k_1$,  then the \emph{conformal fundamental forms}
\begin{gather*}
\alpha_1  =  \sqrt{(k_3 - k_1)(k_2 - k_1)}\omega_1, \qquad
\alpha_2  =  \sqrt{(k_2 - k_1)(k_3 - k_2)}\omega_2, \\
\alpha_3  =  \sqrt{(k_3 - k_2)(k_3 - k_1)}\omega_3
\end{gather*}
are closed, if and only if, the hypersurface $M^3$ is conformally f\/lat. Therefore, when $\alpha_i$ are closed forms, locally there exist $x_1$, $x_2$, $x_3$ such that $\alpha_1 = dx_1$, $\alpha_2 = dx_2$ and $\alpha_3 = dx_3$. By integration, we obtain a special principal coordinate system $x_1$, $x_2$, $x_3$ for a conformally f\/lat hypersurface in~$M^4_K$.

\begin{definition}\label{def.guichard}
A triply orthogonal coordinate system in a Riemannian 3-manifold $(M, g)$
\begin{gather*}
x = (x_1, x_2, x_3) : \ (M, g) \rightarrow \R^3,
\end{gather*}
where the functions  $l_i = \sqrt{g \left(\partial_{x_i}, \partial_{x_i} \right)}$ satisfy the Guichard condition
\begin{gather}
l_1^2 - l_2^2 + l_3^2 = 0, \label{relacaoguichard}
\end{gather}
 is called  a \textit{Guichard net}.
\end{definition}

Since we can deform a conformally f\/lat immersion in a space form into a f\/lat immersion in the light cone, we can consider Guichard nets for f\/lat immersions $f : M^3 \rightarrow L^5$. For such a f\/lat immersion, we express the induced metric $g = \langle  df, df \rangle $, in terms of the Guichard net, as
\begin{gather*}
g = l_1^2 dx_1^2 + l_2^2 dx_2^2 + l_3^2 dx_3^2.
\end{gather*}
Since the metric is f\/lat, the functions $l_i$ must satisfy the {\it Lam\'e's system} \cite[pp.~73--78]{reflame}:
\begin{gather}
	\frac{\partial^2 l_i}{\partial x_j \partial x_k} - \frac{1}{l_j}\frac{\partial l_i}{\partial x_j} \frac{\partial l_j}{\partial x_k} - \frac{1}{l_k}\frac{\partial l_i}{\partial x_k}\frac{\partial l_k}{\partial x_j}  =  0, \nonumber\\
	 \frac{\partial}{\partial x_j} \left( \frac{1}{l_j} \frac{\partial l_i}{\partial x_j} \right) + \frac{\partial}{\partial x_i} \left( \frac{1}{l_i} \frac{\partial l_j}{\partial x_i} \right) + \frac{1}{l_k^2}\frac{\partial l_i}{\partial x_k} \frac{\partial l_j}{\partial x_k}  = 0.
\label{lame}
 \end{gather}
for $i$, $j$, $k$ distinct. Moreover, if $f : M^3 \rightarrow L^5$ is f\/lat, we can consider $M^3$ as a subset of the Euclidean space $\R^3$ and $f$ as an isometric immersion. Then we have a Guichard net on an open subset of $\R^3$, by considering as in Def\/inition~\ref{def.guichard}, $x :  U \subset \R^3 \rightarrow \R^3$, where the functions $l_i$ satisfy the Guichard condition (\ref{relacaoguichard}) and the Lam\'e's system~(\ref{lame}). At this point, one can ask if such a~Guichard net determines a conformally f\/lat hypersurface in a space form, or equivalently, a~f\/lat immersion in~$L^5$.  The answer to this question was given by the following fundamental result due to Hertrich-Jeromin~\cite{Jeromin1}:

\begin{theorem}\label{existence}
For any generic conformally flat hypersurface of a space form $M^4_K$, there exists a~Guichard net $x : U \subset \R^3 \rightarrow \R^3$ on an open set~$U$ of the Euclidean space $\R^3$ $($uniquely determined up to a M{\"o}bius transformation of $\R^3)$.

 Conversely, given a Guichard net $x = (x_1,   x_2,   x_3)  : U \subset \R^3 \rightarrow \R^3$ for the Euclidean space, with $l_i = \sqrt{g(\partial_{x_i}, \partial_{x_j})}$, where $g$ is the canonical flat metric, there exists a generic conformally flat hypersurface in a space form $M^4_K$ $($in this case, M{\"o}bius equivalent Guichard nets are related to conformally equivalent immersions$)$, whose induced metric is given by
\begin{gather}
g = e^{2 P(x)} \big\{ l_1^2 dx_1^2 + l_2^2 dx_2^2 + l_3^2 dx_3^2  \big\}, \label{metricexistence}
\end{gather}
where $P(x)$ is a function depending on~$M^4_K$.
\end{theorem}

The converse is based on the fact that the functions $l_i$ determine the connection forms of a f\/lat immersion $f : M^3 \rightarrow L^5$. In fact, these connection forms satisfy the Maurer--Cartan equations  if,  and only if, the functions $l_i$ satisfy the Guichard condition and the Lam\'e's system.

Therefore, one way of obtaining generic conformally f\/lat hypersurfaces in space forms $M^4_K$ is  f\/inding solutions of  Lam\'e's system, satisfying the Guichard condition. Then the hypersurfaces are constructed by using  Theorem~\ref{existence}. Our objective is to obtain a class of such solutions and to investigate the associated Guichard nets as well as the conformally f\/lat hypersurfaces. We will use the theory of Lie point symmetry groups of dif\/ferential equations, to obtain the symmetry group of Lam\'e's system and their solutions invariant under the action of subgroups of the symmetry group. This is the content of the following sections.

\section{The symmetry group of Lam\'e's system}\label{section3}

In this section, we obtain the Lie point symmetry group of Lam\'e's system. We start with a~brief introduction of symmetry groups of dif\/ferential equations. The reader who is familiar to the theory may skip this introduction.

The theory of Lie point symmetry group is an important tool for the analysis of dif\/ferential equations developed by Lie at the end of the nineteen century \cite{lie}. Roughly speaking, Lie point symmetries of a system of dif\/ferential equations consist of a Lie group of transformations acting on the dependent and independent variables, that transform solutions of the system into solutions.

A standard reference for the theory of symmetry groups of dif\/ferential equations is Olver's book \cite{Olver}, where a clear approach to the subject is given, with theoretical foundations and a~large number of examples and techniques. We will describe here some basic concepts that will be used in this section.

A system $S$ of $n$-th order dif\/ferential equations in $p$ independent and $q$ dependent variables is given as a system of equations
\begin{gather}
\Delta_{r} \big(x, u^{(n)}\big) = 0, \qquad v = 1, \ldots, l, \label{genericsystem}
\end{gather}
involving $x=(x_1,\ldots, x_p)$, $u=(u_1, \ldots, u_q)$ and the derivatives $u^{(n)}$ of $u$ with respect to $x$ up to order $n$.

A \emph{symmetry group} of the system $S$ is a local Lie group of transformations $G$ acting on an open subset $M \subset X \times U$ of the space of independent and dependent variables for the system, with the property that whenever $u = f(x)$ is a solution of $S$, and whenever $gf$ is def\/ined for $g \in G$, then $u = g f(x)$ is also a solution of the system. A vector f\/ield $\textbf{v}$ in the Lie algebra $\mathfrak{g}$ of the group $G$ is called an \emph{infinitesimal generator}.

Consider $\textbf{v}$ as a vector f\/ield on $M \subset X \times U$, with corresponding (local) one-parameter group $\exp(\varepsilon \textbf{v})$, i.e.,
\begin{gather*}
\exp (\varepsilon  \textbf{v}) \equiv \Psi (\varepsilon, x), 
\end{gather*}
where $\Psi$ is the \emph{flow} generated by $\textbf{v}$. In this case, $\textbf{v}$ will be the inf\/initesimal generator of the action.

The symmetry group of a given system of dif\/ferential equation, is obtained by using the \emph{prolongation formula} and the \emph{infinitesimal criterion} that are described as follows. Given a~vector f\/ield on $M \subset X \times U$,
\begin{gather*}
\textbf{v} =   \sum_{i=1}^p \xi^i (x, u) \frac{\partial}{\partial x_i} +  \sum_{\alpha = 1}^q \phi_{\alpha} (x, u) \frac{\partial}{\partial u^{\alpha}},
\end{gather*}
 the \emph{$n$-th prolongation} of $\textbf{v}$ is the vector f\/ield
\begin{gather*}
{\rm pr}^{(n)} \textbf{v} = \textbf{v} +   \sum_{\alpha =1}^q \sum_{J} \phi_{\alpha}^J \big(x, u^{(n)}\big) \frac{\partial}{\partial u^{\alpha}_J}.
\end{gather*}
It is def\/ined on the corresponding \emph{jet space} $M^{(n)} \subset X \times U^{(n)}$, whose coordinates represent the independent variables, the dependent variables and the derivatives of the dependent variables up to order $n$. The second summation is taken over all (unordered) multi-indices $J = (j_1, \ldots, j_k)$, with $1 \leq j_k \leq p$, $1 \leq k \leq n$. The coef\/f\/icient functions $\phi_{\alpha}^J$ of ${\rm pr}^{(n)} \textbf{v}$ are given by the following formula:
\begin{gather*}
\phi_{\alpha}^J \big(x, u^{(n)}\big) = D_J \left(\phi_{\alpha} -   \sum_{i=1}^p \xi^i u^{\alpha}_{J, i}, \right),
\end{gather*}
where $u^{\alpha}_i = \frac{\partial u^{\alpha}}{\partial x_i}$, $u^{\alpha}_{J, i} = \frac{\partial u^{\alpha}_J}{\partial x_i}$ and $D_J$ is given by the total derivatives
\begin{gather*}
D_J = D_{j_1} D_{j_2} \cdots D_{j_k},
\end{gather*}
with
\[
D_i f\big(x, u^{(n)}\big) = \frac{\partial f}{\partial x_i} +   \sum\limits_{\alpha_1}^p   \sum_J u^{\alpha}_{J, i} \frac{\partial f}{\partial u^{\alpha}_J}.
\]

We say that the system (\ref{genericsystem}) is a system of \emph{maximal rank} over $M \subset X \times U$, if the Jacobian matrix
\begin{gather*}
J_{\Delta} \big(x, u^{(n)}\big) = \left(\frac{\partial \Delta_{r}}{\partial x_i}, \frac{\partial \Delta_{r}}{\partial u^{\alpha}_{,\textbf{J}}}\right)
\end{gather*}
has rank $l$, whenever $\Delta_{r} \big(x, u^{(n)}\big) = 0$, where $\textbf{J} = (j_1, \ldots, j_k)$ is a multi-index that denotes the partial derivatives of $u^{\alpha}$.

Suppose that (\ref{genericsystem}) is a system of maximal rank. Then the set of all vectors f\/ields $\textbf{v}$ on $M$ such that
\begin{gather}
{\rm pr}^{(n)}  \textbf{v}  \big[\Delta_{r} \big(x, u^{(n)}\big)\big] = 0, \qquad  r=1, \ldots, l, \qquad {\rm whenever} \quad \Delta_{r} \big(x, u^{(n)}\big) = 0, \label{criterion}
\end{gather}
is a Lie algebra of inf\/initesimal generators of a symmetry group for the system. It is shown in~\cite{Olver} that the inf\/initesimal criterion~(\ref{criterion}) is in fact both a necessary and suf\/f\/icient condition for a group~$G$ to be a symmetry group. Hence, all the connected symmetry groups can be determined by considering this criterion.

 Since the prolongation formula is given in terms of $\xi^i$ and $\phi_{\alpha}$ and  the partial derivatives  with respect to both~$x$ and~$u$, the inf\/initesimal criterion provides a system of partial dif\/ferential equations for the coef\/f\/icients $\xi^i$ and $\phi_{\alpha}$ of~$\textbf{v}$, called the \emph{determining equations}.  By solving these equations, we obtain the vector f\/ield $\textbf{v}$ that determines a Lie algebra $\mathfrak{g}$. The symmetry group~$G$ is obtained by exponentiating the Lie algebra.

\subsection{Obtaining the symmetry group of Lam\'e's system}\label{section3.1}

From now on, we consider the following notation for  derivatives of a function $f= f(x_1, \ldots, x_n)$
\begin{gather*}
f_{,x_i}:=\frac{\partial f}{\partial x_i} \qquad  {\rm and} \qquad f_{,x_ix_j} := \frac{\partial^2 f}{\partial x_i \partial x_j}.
\end{gather*}
With this notation, Lam\'e's system (\ref{lame}) is given by
\begin{gather}
	l_{i,x_jx_k} - \frac{l_{i,x_j} l_{j,x_k}}{l_j} - \frac{l_{i,x_k} l_{k,x_j}}{l_k} = 0, \label{notalame1} \\
	\left( \frac{l_{i,x_j}}{l_j} \right)_{,x_j} + \left( \frac{l_{j,x_i}}{l_i} \right)_{,x_i} + \frac{l_{i,x_k} l_{j,x_k}}{l_k^2} = 0, \label{notalame2}
\end{gather}
where $i$, $j$ and $k$ are distinct indices in the set $\left\{ 1,2,3\right\}$. We will also consider the following notation,
\begin{gather}
	\varepsilon_s = \left\{
	   \begin{array}{@{}rl}
            1  & {\rm if}  \ s=1 \ {\rm or} \ s=3, \\
            -1 & {\rm if}  \ s=2.
     \end{array}
	\right.
	\label{sinalepsilon}
\end{gather}
We can now rewrite Guichard condition as
\begin{gather*}
	\varepsilon_i l_i^2 + \varepsilon_j l_j^2 + \varepsilon_k l_k^2 = 0.
\end{gather*}

Next, we introduce auxiliary functions in order to reduce the system of second-order dif\/fe\-ren\-tial equations  (\ref{notalame1}) and (\ref{notalame2}), into a f\/irst order one. Consider the functions $h_{ij}$, with $i \neq j$, given by
\begin{gather*}
	l_{i,x_j} - h_{ij}{l_j} =0.
\end{gather*}
With these functions,  we rewrite (\ref{notalame1}) and (\ref{notalame2}) as
\begin{gather*}
	h_{ij,x_k} - h_{ik} h_{kj}  =  0,  \qquad
	h_{ij,x_j} + h_{ji,x_i} + h_{ik} h_{jk}  =  0.
\end{gather*}
for $i$, $j$, $k$ distinct. Since the functions $l_1$, $l_2$ and $l_3$ satisfy  Guichard condition, there are other relations involving the derivatives of~$l_i$ and $h_{ij}$. Taking the derivative of Guichard condition  with respect to~$x_i$, we have
\begin{gather*}
	\varepsilon_i l_{i,x_i} + \varepsilon_j h_{ji} l_j + \varepsilon_k h_{ki} l_k = 0,
\end{gather*}
for $i$, $j$, $k$ distinct. The derivatives of  the above equation with respect to $x_j$ leads to
\begin{gather*}
	\varepsilon_i h_{ij,x_i} + \varepsilon_j h_{ji,x_j} + \varepsilon_k h_{ki} h_{kj} = 0.
\end{gather*}

Therefore, we summarize the last six equations in the following system of f\/irst-order partial dif\/ferential equations, equivalent to Lam\'e's system, that we call \textit{Lam\'e's system of first order}
\begin{gather}
\varepsilon_i l_i^2 + \varepsilon_j l_j^2 + \varepsilon_k l_k^2 = 0, \label{A} \\
l_{i,x_j} - h_{ij}l_j = 0, \label{B} \\
\varepsilon_i l_{i,x_i} + \varepsilon_j h_{ji}l_j + \varepsilon_k h_{ki}l_k = 0, \label{C} \\
h_{ij,x_k} - h_{ik} h_{kj} = 0, \label{D} \\
h_{ij,x_j} + h_{ji,x_i} + h_{ik}h_{jk} = 0, \label{E} \\
\varepsilon_i h_{ij,x_i} + \varepsilon_j h_{ji,x_j} + \varepsilon_k h_{ki}h_{kj} = 0. \label{F}
\end{gather}
By considering $x = (x_1, x_2, x_3)$, $l = (l_1, l_2, l_3)$ and $h$ the of\/f-diagonal $3 \times 3$ matrix given by $h_{ij}$ in our next two results, we obtain the Lie algebra of the inf\/initesimal generators and the symmetry group of Lam\'e's system of f\/irst order.
\begin{theorem} \label{teo.generator}
Let  $V$ be the infinitesimal generator of the symmetry group of Lam\'e's system of first order \eqref{A}--\eqref{F}, given by
\begin{gather}
V =  \sum_{i=1}^3 \xi^i(x,l,h) \frac{\partial}{\partial x_i} +  \sum_{i=1}^3 \eta^i(x,l,h) \frac{\partial}{\partial l_i} +  \sum_{i,j=1, \, i \neq j}^3 \phi^{ij}(x,l,h) \frac{\partial}{\partial h_{ij}}. \label{infgeneratorlame}
\end{gather}
Then the functions $\xi^i$, $\eta^i$ and $\phi^{ij}$ are given by
\begin{gather*}
\xi^i = a x_i + a_i, \qquad
\eta^i = c l_i, \qquad
\phi^{ij} = - a h_{ij},
\end{gather*}
 where  $a, c,a_i \in \R$.
\end{theorem}

The proof of Theorem~\ref{teo.generator} is very long and technical. It consists of obtaining the functions~$\xi^i$,~$\eta^i$ and~$\phi^{ij}$ by solving the determining equations which are obtained as follows. We apply the f\/irst prolongation of~$V$ to each equation (\ref{A})--(\ref{F}) and we eliminate the functional dependence of the derivatives of $h$ and $l$ caused by the system. Then we equate to zero the coef\/f\/icients of the remaining unconstrained partial derivatives. The complete proof with, all the details, is given in Appendix~\ref{appendix}.

As a consequence of Theorem~\ref{teo.generator}, by exponentiating $V$, we obtain the symmetry group of Lam\'e's system. Observe that the functions~$\phi^{ij}$ do not depend on $x$ and $l$ (see \cite{olver2} for symmetry group of equivalent systems):
\begin{corollary}\label{corollary1}
The symmetry group of Lam\'e's system \eqref{A}--\eqref{F} is given by the following transformations:
\begin{enumerate}\itemsep=0pt
\item[$1)$] translations in the independent variables: $\tilde{x_i} = x_i + v_i$;
\item[$2)$] dilations in the independent variables: $\tilde{x_i} = \lambda x_i$;
\item[$3)$] dilations in the dependent variables: $\tilde{l_i} = \rho l_i$;
\end{enumerate}
where $v_i\in \R$ and $\lambda, \rho \in \R\setminus \{0\}$.
\end{corollary}

\subsection{Group invariant solutions}\label{section3.2}

The knowledge of all the inf\/initesimal generators $\textbf{v}$ of the symmetry group of a system of dif\/ferential equations, allows one to reduce the system to another one with a reduced number of variables. Specif\/ically, if the system has $p$ independent variables and an $s$-dimensional symmetry subgroup is considered, then the reduced system for the  solutions invariant under this subgroup will depend on $p-s$ variables (see Olver~\cite{Olver} for details). Finding all the $s$-dimensional symmetry subgroups is equivalent to f\/inding all the $s$-dimensional subalgebras of the Lie algebra of inf\/initesimal symmetries $\textbf{v}$. For the remainder of this paper, we will consider the 2-dimensional subgroups of the symmetry group of Lam\'e's system. The f\/irst one will be the  translation subgroup and the second one will be the subgroup involving translations and the dilations. The 1-dimensional subgroup given just by dilations and the solutions invariant under this subgroup are being investigated. We will report on our investigation in another paper. We observe that the symmetry subgroup of dilations in the dependent variables (Corollary~\ref{corollary1}(3)) is irrelevant for the geometric study of conformally f\/lat hypersurfaces due to~(\ref{metricexistence}).

We start with the 2-dimensional subgroup of translations. The basic invariant of this group is given by
\begin{gather}
\xi = \alpha_1 x_1 + \alpha_2 x_2 + \alpha_3 x_3, \label{invarianteconjunto}
\end{gather}
where $(\alpha_1, \alpha_2, \alpha_3)$ is a non zero  vector. We will consider solutions  $l_i$ such that
\begin{gather*}
l_i(x_1, x_2, x_3) = l_i (\xi), \qquad 1 \leq i \leq 3,  
\end{gather*}
where $\xi$ is given by (\ref{invarianteconjunto}). For such solutions,  Lam\'e's system reduces to a system of ODEs. We start with two lemmas:

\begin{lemma} \label{cieli}
Let $l_s (\xi)$, $s=1, 2,3,$ where $\xi = \sum\limits_{s=1}^3 \alpha_s x_s$, be a solution of Lam\'e's system \eqref{A}--\eqref{F}. Let $i,  k \in \left\{ 1,  2,  3 \right\}$ be two fixed and distinct indices such that $\alpha_i = \alpha_k = 0$. Then $l_i$ or $l_k$ is constant.
\end{lemma}

\begin{proof}
Since $\alpha_i = \alpha_k = 0$, it follows from (\ref{B}) that equation (\ref{E}) reduces to
\begin{gather*}
\alpha_j^2 \left[ \frac{l_{i,\xi}}{l_j} \right]_{,\xi} = 0,
\end{gather*}
which implies $l_{i,\xi} = c_i l_j$, where $c_i \in \R$. Similarly, interchanging $i$ with $k$, we obtain  $l_{k,\xi} = c_k l_j$. Finally, interchanging~$k$ with  $j$, we get
\begin{gather*}
\alpha_j^2 \frac{l_{i,\xi}l_{k,\xi}}{l_j^2} = \alpha_j^2 c_i c_k = 0.
\end{gather*}
Therefore, we conclude that $l_i$ or $l_k$ is constant.
\end{proof}

\begin{lemma}\label{ljalphaj}
Let $l_s (\xi)$, $s=1, 2,3,$ where $\xi =  \sum\limits_{s=1}^3 \alpha_s x_s$, be a solution of Lam\'e's system \eqref{A}--\eqref{F}. If there exists a unique $j \in \left\{ 1,  2,  3 \right\}$ such that $l_j$ is a non zero constant, then  $\alpha_j = 0$.
\end{lemma}

\begin{proof}
Interchanging the indices in $(\ref{D})$, we obtain the following two equations
\begin{gather}
	    \alpha_j \alpha_k \left( l_{i,\xi \xi} - \frac{l_{i,\xi}l_{k, \xi}}{l_k} \right) = 0  \label{d1}, \\
	    \alpha_j \alpha_i \left( l_{k,\xi \xi} - \frac{l_{k,\xi}l_{i, \xi}}{l_i} \right) = 0  \label{d2},
\end{gather}
and an identity.

Similarly, it follows from $(\ref{E})$ that
\begin{gather}	
	\alpha_j^2 l_{i, \xi \xi} = 0 \label{e1}, \\
	\alpha_j^2 l_{k, \xi \xi} = 0 \label{e2}, \\
	\alpha_k^2 \left( \frac{l_{i,\xi}}{l_k} \right)_{,\xi}
+ \alpha_i^2 \left( \frac{l_{k,\xi}}{l_i} \right)_{,\xi} + \alpha_j^2 \frac{l_{i,\xi}l_{k,\xi}}{l_j} = 0. \label{e3}
\end{gather}

Suppose, by contradiction,  that $\alpha_j \neq 0$. It follows from (\ref{e1}) and (\ref{e2}) that $l_{i, \xi} = c_i$ and $l_{k, \xi} = c_k$, where $c_i \neq 0$ and $c_k \neq 0$, since by hypothesis, $l_i$ and $l_k$ are non constants.  Then, it follows from (\ref{d1}) and (\ref{d2}) that $\alpha_i = \alpha_k = 0$. From (\ref{e3}),  we obtain $ \alpha_j^2 c_i c_k = 0,$ which is a~contradiction.
\end{proof}

The following theorem gives the solutions of Lam\'e's
system, satisfying Guichard condition, which are invariant under the action of the translation group, whenever none of the functions $l_i$ is constant.

\begin{theorem} \label{teo.nonconstant.trans}
Let $l_s (\xi)$, $s=1, 2,3,$ where $\xi = \sum\limits_{s=1}^3 \alpha_s x_s$, be a solution of Lam\'e's system \eqref{A}--\eqref{F}, such that $l_{s}$ is not constant for all $s$. Then there exist $c_s \in \R \setminus \left\{ 0 \right\}$, such that,
    \begin{gather}
    l_{i, \xi} = c_i l_k l_j, \quad i,  j,   k \ \ \text{distinct} \label{derivadalame}, \\
	c_1 - c_2 + c_3 = 0, \label{apenasce} \\
	\alpha_1^2 c_2 c_3 + \alpha_2^2 c_1 c_3 + \alpha_3^2 c_1 c_2 = 0. \label{alface}
	\end{gather}
Moreover, the functions $l_i(\xi)$ are given by
	\begin{gather}
	l_{1,\xi}^2 = c_2(c_2 - c_1) \left( l_1^2 - \frac{\lambda}{c_2} \right) \left( l_1^2 - \frac{\lambda}{c_2 - c_1}  \right),  \label{ode.jacobi} \\
	l_2^2 = \frac{c_2}{c_1} \left( l_1^2 - \frac{\lambda}{c_2} \right), \label{l2.jacobi} \\
	l_3^2 = \frac{c_2 - c_1}{c_1} \left(l_1^2 - \frac{\lambda}{c_2 - c_1} \right), \label{l3.jacobi}
	\end{gather}
where $\lambda\in \R$.
\end{theorem}
\begin{proof}
 By hypothesis, we are considering non constant solutions. Then, it follows from Lem\-ma~\ref{cieli}, that~$\alpha_s \neq 0$ for at least two distinct indices. Suppose that $\alpha_j$ and $\alpha_k$ non zero. From~(\ref{B}) and~(\ref{D}) we obtain
\begin{gather*}
\alpha_j \alpha_k \left\{ \left[ \frac{l_{i,\xi}}{l_j} \right]_{,\xi} - \frac{l_{i,\xi}}{l_j} \frac{l_{k,\xi}}{l_k}  \right\} = 0,
\end{gather*}
which implies
\begin{gather*}
	\left[\frac{l_{i,\xi}}{l_j} \right]_{,\xi}\left[\frac{l_{i,\xi}}{l_j} \right]^{-1} = \frac{l_{k,\xi}}{l_k} \label{intdolog}.
\end{gather*}
Integrating this equation, we obtain $l_{i,\xi} = c_i l_k l_j$, where $c_i \neq 0$.

If $\alpha_i \neq 0$, analogously considering the non zero pairs $(\alpha_i ,   \alpha_j)$ and $(\alpha_i,   \alpha_k)$, we conclude that $l_{k, \xi} = c_k l_i l_j$ and $l_{j, \xi} = c_j l_i l_k$.  If  $\alpha_i = 0$, then from equation (\ref{E}) we have
\begin{gather*}
\left[\frac{l_{i,x_j}}{l_j} \right]_{,x_j} + \frac{l_{i,x_k}}{l_k} \frac{l_{j,x_k}}{l_k} = \alpha_j^2 c_i l_{k,\xi} + \alpha_k^2 c_i l_j \frac{l_{j,\xi}}{l_k} = 0.
\end{gather*}
Since $c_i \neq 0$, we integrate the above expression to obtain
\begin{gather*}
\alpha_j^2 l_k^2 + \alpha_k^2 l_j^2 = \lambda_{jk},
\end{gather*}
where $\lambda_{jk}$ is a constant. This equation and Guichard condition (\ref{A}) lead to
\begin{gather}
l_j^2  =  \frac{\alpha_j^2}{\alpha_k^2} \left(\frac{\lambda_{jk}}{\alpha_j^2} - l_k^2 \right),  \qquad
l_i^2  =  \frac{\ei}{\alpha_k^2} \big[ l_k^2 \big(\ej \alpha_j^2 - \ek \alpha_k^2 \big) - \ej \lambda_{jk}  \big].
	\label{ljfunclk}
\end{gather}
Taking the derivative of the last equation with respect to $\xi$, we have
\begin{gather*}
l_i \left( c_i l_k l_j \right) = \frac{\ei}{\alpha_k^2} \big[ l_k l_{k,\xi} \big(\ej \alpha_j^2 - \ek \alpha_k^2 \big) \big].
\end{gather*}
If $\ej \alpha_j^2 - \ek \alpha_k^2 \neq 0$, we conclude that
\begin{gather*}
 l_{k,\xi} = \frac{c_i \alpha_k^2}{\ej \alpha_j^2 - \ek \alpha_k^2} l_i l_j = c_k l_i l_j.
\end{gather*}
Applying this expression into the derivative of the f\/irst equation in (\ref{ljfunclk}) with respect to  $\xi$ we obtain
\begin{gather*}
l_j l_{j, \xi} = - \frac{\alpha_j^2}{\alpha_k^2}l_k l_{k, \xi} =  - \frac{\alpha_j^2}{\alpha_k^2}l_k \left(c_k l_i l_j \right),
\end{gather*}
consequently, $ l_{j, \xi} = c_j l_i l_k$.

Next, we will show that $\ej \alpha_j^2 - \ek \alpha_k^2 \neq 0$ to conclude the proof of~(\ref{derivadalame}). Suppose by contradiction that $\ej \alpha_j^2 - \ek \alpha_k^2 = 0$, then the f\/irst equation of~(\ref{ljfunclk}) can be written as $ \ej l_j^2 + \ek l_k^2 = \frac{\ek \lambda_{jk}}{\alpha_j^2} $. Then Guichard condition now implies  that $l_i$ is constant, which is a contradiction. The relations between the constants~(\ref{apenasce}) and~(\ref{alface}) follow from a straightforward computation using equations  (\ref{C}) and (\ref{E}), respectively.

In order to complete the proof of the theorem, we start with
\begin{gather}
l_{1, \xi} = c_1 l_2 l_3, \label{l1qsi}  \\
l_{2, \xi} = c_2 l_1 l_3, \label{l2qsi}   \\
l_{3, \xi} = c_3 l_1 l_2. \label{l3qsi}
\end{gather}
Multiplying (\ref{l2qsi}) by  $l_2$ and integrating we have
\begin{gather}
	l_2^2 = \frac{c_2}{c_1} \left( l_1^2 - \frac{\lambda}{c_2} \right),
	\label{l2funcl1}
\end{gather}
where $\lambda$ is a constant. Therefore, it follows from (\ref{l2funcl1}) and Guichard condition that
\begin{gather}
	l_3^2 = \frac{c_2 - c_1}{c_1} \left( l_1^2 - \frac{\lambda}{c_2 - c_1} \right).
	\label{l3funcl1}
\end{gather}
Using (\ref{l1qsi}), (\ref{l2funcl1}) and (\ref{l3funcl1}), we conclude that
\begin{gather*}
	{l_{1, \xi}}^2 = c_1^2 \left[ \frac{c_2}{c_1} \left( l_1^2 - \frac{\lambda}{c_2} \right) \right] \left[ \frac{c_2 - c_1}{c_1} \left( l_1^2 - \frac{\lambda}{c_2 - c_1} \right) \right] \\
\hphantom{{l_{1, \xi}}^2 }{}
	= c_2 \left(c_2 - c_1 \right) \left( l_1^2 - \frac{\lambda}{c_2} \right) \left( l_1^2 - \frac{\lambda}{c_2 - c_1} \right).\tag*{\qed}
\end{gather*}
\renewcommand{\qed}{}
\end{proof}

In our next theorem,  we consider the solutions $l_i (\xi)$ when one of the functions  $l_i$ is constant.
\begin{theorem} \label{uma.constante}
Let $l_s (\xi)$, $s=1,  2, 3,$ where $\xi =  \sum\limits_{s=1}^3 \alpha_s x_s$, be a solution of Lam\'e's system \eqref{A}--\eqref{F}. Suppose that only one of the functions $l_s$ is constant. Then one of the following occur:
\begin{enumerate}\itemsep=0pt
\item[$a)$] $l_1 = \lambda_1$, $l_2 = \lambda_1 \cosh (b\xi + \xi_0)$, $l_3 = \lambda_1 \sinh (b\xi + \xi_0)$, where $\xi = \alpha_2 x_2 + \alpha_3 x_3$, $\alpha_2^2 + \alpha_3^2 \neq 0$ and $b,  \xi_0\in \R$ ;
\item[$b)$] $ l_2 = \lambda_2$, $l_1 = \lambda_2 \cos \varphi(\xi)$, $l_3 = \lambda_2 \sin \varphi(\xi) $, where $\xi = \alpha_1 x_1 + \alpha_3 x_3$, $\alpha_1^2 + \alpha_3^2 \neq 0$ and $\varphi$ is one  of the following:
   \begin{enumerate}\itemsep=0pt
      \item[$b.1)$] $\varphi(\xi) = b \xi + \xi_0$, if $\alpha_1^2 \neq \alpha_3^2$, where $\xi_0, b\in \R$;
      \item[$b.2)$] $\varphi$ is any function of $\xi$, if $\alpha_1^2 = \alpha_3^2$;
   \end{enumerate}
\item[$c)$]  $l_3 = \lambda_3$, $l_2 = \lambda_3 \cosh (b\xi + \xi_0)$, $l_1 = \lambda_3 \sinh (b\xi + \xi_0)$, where $\xi = \alpha_1 x_1 + \alpha_2 x_2$, $\alpha_1^2 + \alpha_2^2 \neq 0$ and   $b,\xi_0\in\R $.
\end{enumerate}
\end{theorem}

\begin{proof}
We will consider each case separately:

$a)$ If $l_1 = \lambda_1$, then it follows from Lemma \ref{ljalphaj} that we must have $\xi = \alpha_2 x_2 + \alpha_3 x_3$. Now  Guichard condition implies that $l_2 = \lambda_1 \cosh \varphi(\xi)$ and $l_3 = \lambda_1 \sinh \varphi(\xi)$. In order to determi\-ne~$\varphi$, we use~(\ref{E}) with the following indices
\begin{gather*}
h_{23,x_3}+ h_{32,x_2} + h_{21} h_{31} = 0,
\end{gather*}
to obtain
\begin{gather*}
\alpha_3^2 \left( \frac{\lambda_1 \varphi_{,\xi} \sinh \varphi}{\lambda_1 \sinh \varphi} \right)_{,\xi} + \alpha_2^2 \left( \frac{\lambda_1 \varphi_{,\xi} \cosh \varphi}{\lambda_1 \cosh \varphi} \right)_{,\xi} = 0.
\end{gather*}
Since $l_2$ and $l_3$ are not constant, we have $\alpha_2^2 + \alpha_3^2 \neq 0$, which implies $\varphi_{,\xi \xi} = 0.$ Consequently, $\varphi (\xi) = b\xi + \xi_0$.

$b)$ If $l_2 = \lambda_2$, it follows from Lemma \ref{ljalphaj} that $\xi = \alpha_1 x_1 + \alpha_3 x_3$. Then Guichard condition implies that  $l_1 = \lambda_2 \cos \varphi(\xi)$ and $l_3 = \lambda_2 \sin \varphi(\xi)$. As in the case $a)$, from equation  (\ref{E}) we get $\left(\alpha_1^2 - \alpha_3^2 \right) \varphi_{,\xi \xi} = 0$. Since $l_1$ and $l_3$ are non constant, we have $\alpha_1^2 + \alpha_3^2 \neq 0$. Then we have two cases to consider:

$b.1)$ If $\alpha_1^2 \neq \alpha_3^2$, then  $\varphi(\xi) = b\xi + \xi_0$;

$b.2)$ If $\alpha_1^2 = \alpha_3^2$, then $\varphi$ can be any function of $\xi$.

$c)$ The proof is the same as in $a)$.
\end{proof}

Next, we consider the  solutions invariant under the 2-dimensional subgroup involving translations and dilations. In this case, the basic invariant is given by
\begin{gather}
\eta = \frac{a_1 x_1 + a_2 x_2 + a_3 x_3 }{b_1 x_1 + b_2 x_2 + b_3 x_3},  \label{etainvariant}
\end{gather}
where the vectors $(a_1 ,  a_2 ,  a_3)$ and $(b_1 ,  b_2 ,  b_3)$ are linearly independent. If $f=f(\eta)$ is a function depending on $\eta$, then
\begin{gather*}
f_{,x_i} = f_{,\eta} \eta_{x_i} = \frac{a_i - b_i \eta}{b_1 x_1 + b_2 x_2 + b_3 x_3}f_{,\eta}.
\end{gather*}
In order to simplify the computations, we will use the following notation:
\begin{gather}
N_i := a_i - b_i \eta  \qquad {\rm and} \qquad  \beta = b_1 x_1 + b_2 x_2 + b_3 x_3. \label{notadilatacao}
\end{gather}
 Then we have $\eta_{,x_i} = \frac{N_i}{\beta}$.

In order to obtain the solutions of Lam\'e's system
$l_i(\eta)$, which depend on $\eta$, we will need some lemmas.

\begin{lemma} \label{lema.derivada.dilata}
Let $l_1 (\eta)$, $l_2 (\eta)$, $l_3 (\eta)$, where $\eta$ is given by \eqref{etainvariant}, be a solution of Lam\'e's system \eqref{A}--\eqref{F}. Suppose that  for a fixed pair $j,  k \in \left\{1,  2,  3 \right\}$, $j \neq k$, $(a_j,   b_j) \neq (0, 0)$ and $(a_k,   b_k) \neq (0,  0)$. Then there exists $c_i\in \R $ such that
\begin{gather}
l_{i,\eta} = c_i \frac{l_k l_j}{N_k N_j}, \qquad i \neq j,  k, \label{derivada.dilata}
\end{gather}
where $N_k$ is given by \eqref{notadilatacao}.
\end{lemma}

\begin{proof} From (\ref{B}), we have that $h_{ij} = \frac{l_{i,\eta} N_j}{l_j \beta}$. Then, equation (\ref{D}) can be written as
\begin{gather*}
\left[ \frac{l_{i,\eta} N_k N_j}{l_j} \right]_{\eta} - \frac{l_{i,\eta} N_k N_j}{l_j} \frac{l_{k,\eta}}{l_k} = 0,
\end{gather*}
which implies
\begin{gather*}
\left( \frac{l_{i,\eta}N_k N_j}{l_k l_j} \right)_{,\eta} = 0.
\end{gather*}
Since $(a_j ,  b_j) \neq (0 , 0)$ and $(a_k ,  b_k) \neq (0 , 0)$, we have that $N_j \neq 0$, $N_k \neq 0$ and the equation~(\ref{derivada.dilata}) holds.
\end{proof}

\begin{lemma}\label{li.constant.dilata}
Let $l_1 (\eta),  l_2 (\eta),  l_3 (\eta)$, where $\eta$ is given by \eqref{etainvariant}, be a solution of Lam\'e's system \eqref{A}--\eqref{F}. If $(a_i, b_i) = (0,  0)$, for some $i \in \left\{1,  2,  3 \right\}$, then~$l_i$ is constant.
\end{lemma}

\begin{proof}
Since the vectors $(a_1,  a_2,  a_3)$ and $(b_1,  b_2,  b_3)$ are linearly independent, if  $(a_i,   b_i) = (0,  0)$ we must have $(a_j,   b_j) \neq (0,  0)$ and $(a_k,   b_k) \neq (0,  0)$ for $i,  j,  k$ distinct and we can use Lemma~\ref{lema.derivada.dilata}. By considering equation~(\ref{E}),  we have
\begin{gather*}
\left( \frac{c_i l_k}{\beta N_k} \right)_{,\eta} + \left( \frac{c_i l_k}{\beta N_k} \right) \frac{l_{j,\eta N_k}}{\beta l_k} = 0,
\end{gather*}
which implies
\begin{gather}
c_i \left[ \frac{l_{k,\eta}N_j}{N_k} - \frac{l_k}{N_k^2} \left( N_k \beta \right)_{,x_j} + \frac{l_j l_{j,\eta}N_k}{l_k N_j} \right] = 0. \label{dilataaibi1}
\end{gather}
By interchanging $j$ with $k$,  we have analogously
\begin{gather}
c_i \left[ \frac{l_{j,\eta}N_k}{N_j} - \frac{l_j}{N_j^2} \left( N_j \beta \right)_{,x_k} + \frac{l_k l_{k,\eta}N_j}{l_j N_k} \right] = 0. \label{dilataaibi2}
\end{gather}
Suppose by contradiction that $c_i \neq 0$. Then, it follows from  (\ref{dilataaibi1}) and (\ref{dilataaibi2}) that
\begin{gather*}
\frac{l_k^2}{N_k^2} \left(N_k \beta \right)_{x_j} = \frac{l_j^2}{N_j^2} \left(N_j \beta \right)_{x_k}. 
\end{gather*}
If $a_i = b_i = 0$, we must have
\begin{gather*}
(a_k b_j - b_k a_j) \left( \frac{l_k^2}{N_k^2} + \frac{l_j^2}{N_j^2} \right) = 0,
\end{gather*}
which is a contradiction since $(a_k b_j - b_k a_j) \neq 0$. Therefore $c_i = 0$ and $l_i$ is constant.
\end{proof}

\begin{lemma} \label{liwithaibi}
Let $l_1 (\eta)$, $ l_2 (\eta)$, $l_3 (\eta)$, with
$\eta$ given by \eqref{etainvariant}, be a solution of Lam\'e's system \eqref{A}--\eqref{F}. If there exists a unique function $l_i$ which is a non zero constant, then $(a_i,  b_i) = (0,  0)$.
\end{lemma}

\begin{proof} Suppose by contradiction that $(a_i,  b_i) \neq (0,  0)$. Since $l_j$ and $l_k$ are not constant, for~$i$,~$j$,~$k$ distinct, it follows from Lemma~\ref{li.constant.dilata}, that we must have $(a_j,  b_j) \neq (0,  0)$ and $(a_k,  b_k) \neq (0,  0)$. Then, Lemma~\ref{lema.derivada.dilata} implies that there are constants $c_i$, $c_j$ and $c_k$ such that
\begin{gather*}
l_{i, \eta} = c_i \frac{l_j l_k}{N_j N_k}   , \qquad l_{j, \eta} = c_j \frac{l_k l_i}{N_k N_i} \qquad {\rm and} \qquad l_{k, \eta} = c_k \frac{l_k l_i}{N_k N_i}.
\end{gather*}
Using equation (\ref{E}) and interchanging the indices we have
\begin{gather}
c_k \frac{l_i l_j}{N_j} - \frac{(a_k b_i - b_k a_i)l_k}{N_k} = 0, \label{dilata.ea} \\
c_j c_k \frac{l_j l_k}{N_j N_k} - \frac{l_i}{N_i^2} \left[c_j (a_i b_k - a_k b_i) + c_k (a_i b_j - b_i a_j) \right] = 0, \label{dilata.eb} \\
c_j \frac{l_i l_k}{N_k} - \frac{(a_j b_i - b_j a_i)l_j}{N_j} = 0. \label{dilata.ec}
\end{gather}
Multiplying equation (\ref{dilata.ea}) by $c_j \frac{N_k}{l_k}$, (\ref{dilata.eb}) by $ \frac{N_i^2}{l_i}$ and  (\ref{dilata.ec}) by $c_k \frac{N_j}{l_j}$, the sum will reduce to
\begin{gather*}
c_j c_k \big[ (l_i l_j N_k)^2 + (l_i l_k N_j)^2 + (l_j l_k N_i)^2 \big] = 0,
\end{gather*}
which is a contradiction. Then, we must have $(a_i,   b_i) = (0,   0)$ and the lemma is proved.
\end{proof}

\begin{remark}  We observe that when all pairs $(a_s, b_s)$ are dif\/ferent from zero, then the proof of Lemma~\ref{liwithaibi} shows that the solution $l_i(\eta)$ of Lam\'e's system is constant.
\end{remark}

We will now obtain the solutions $l_s (\eta)$, when one pair $(a_s,   b_s) = (0,   0)$.

\begin{theorem} \label{uma.constante.dilata}
Let $l_i(\eta)$, with $\eta$ given by \eqref{etainvariant}, be a solution of Lam\'e's system invariant under the $2$-dimensional subgroup involving translation and dilations. Suppose that one of the pairs $(a_s, b_s)=(0,0)$. Then one of the following occur:
\begin{enumerate}\itemsep=0pt
\item[$a)$] If $(a_1,   b_1) = (0,   0)$ then  $l_1 = \lambda_1$, $l_2 = \lambda_1 \cosh \varphi(\eta)$, $l_3 = \lambda_1 \sinh \varphi(\eta)$,  where $\eta = \frac{a_2 x_2 + a_3 x_3 }{b_2 x_2 + b_3 x_3}$ and $\varphi$ is given by
\begin{gather}
\varphi(\eta) = \frac{C_0}{a_2 b_3 - a_3 b_2} \arctan \left[ \frac{b_2^2 + b_3^2}{a_3 b_2 - a_2 b_3} \left(\eta - \frac{a_2 b_2 + a_3 b_3}{b_2^2 + b_3^2} \right) \right] + C_1,  \label{eta.dilata1}
\end{gather}
where $C_0,  C_1\in\R$.
\item[$b)$] If $(a_2,   b_2) = (0,   0)$ then  $ l_2 = \lambda_2$, $l_1 = \lambda_2 \cos \varphi(\eta)$, $l_3 = \lambda_2 \sin \varphi(\eta) $, where $\eta = \frac{a_1 x_1 + a_3 x_3 }{b_1 x_1 + b_3 x_3}$ and~$\varphi$ is given as follows:
   \begin{enumerate}\itemsep=0pt
      \item[$b.1)$] if $b_1 = b_3 = b$, then
\begin{gather}
      \varphi(\eta) = \frac{D_0}{2b(a_3 - a_1)} \log \left( 2b \eta - a_1 - a_3 \right) + D_1, \label{eta.dilata21}
\end{gather}
where $D_0,D_1\in\R$;
      \item[$b.2)$] if $b_1 \neq b_3$, then
\begin{gather}
\varphi (\eta) = \frac{D_2}{2(a_1b_3 - a_3b_1)} \log \left[ \frac{(b_3 + b_1)\eta - (a_3 + a_1)}{(b_3 - b_1)\eta - (a_3 - a_1)} \right] + D_3, \label{eta.dilata22}
\end{gather}
where  $D_2,D_3\in \R$.
   \end{enumerate}
\item[$c)$]  If $(a_3,   b_3) = (0,   0)$, then $l_3 = \lambda_3$, $l_2 = \lambda_3 \cosh\varphi(\eta)$, $l_1 = \lambda_3 \sinh\varphi(\eta)$, with $\eta = \frac{a_1 x_1 + a_2 x_2 }{b_1 x_1 + b_2 x_2}$ and $\varphi$ is given by
\begin{gather*}
\varphi(\eta) = \frac{E_0}{a_2 b_1 - a_1 b_2} \arctan \left[ \frac{b_2^2 + b_1^2}{a_2 b_1 - a_1 b_2} \left(\eta - \frac{a_2 b_2 + a_1 b_1}{b_2^2 + b_1^2} \right) \right] + E_1, 
\end{gather*}
where $E_0, E_1\in \R$.
\end{enumerate}
\end{theorem}

\begin{proof}
$a)$ If $(a_1,   b_1) = (0,   0)$ then Lemma \ref{li.constant.dilata} implies that  $l_1 = \lambda_1$ and Guichard condition implies that $l_2 = \lambda_1 \cosh \varphi (\eta)$ and $l_3 = \lambda_1 \sinh \varphi(\eta)$.
In order to f\/ind $\varphi$, we use equation (\ref{E}) with the following indices
\begin{gather*}
h_{32,x_2} + h_{23,x_3} + h_{31} h_{21} = 0.
\end{gather*}
Since  $h_{32} = \frac{\varphi_{,\eta} N_2}{\beta}$,  $h_{23} = \frac{\varphi_{,\eta} N_3}{\beta}$ and $h_{31}=h_{21} = 0$, we rewrite the equation above as
\begin{gather*}
\left( \frac{\varphi_{,\eta} N_2}{\beta} \right)_{,x_2} + \left( \frac{\varphi_{,\eta} N_3}{\beta} \right)_{,x_3} = 0.
\end{gather*}
By substituting the derivatives,  we have the following ODE
\begin{gather*}
\varphi_{,\eta \eta} \big[ (N_2)^2 + (N_3)^2 \big] - 2 \varphi_{,\eta} (b_2 N_2 + b_3 N_3) = 0,
\end{gather*}
whose solution is exactly~(\ref{eta.dilata1}).

$b)$ If $(a_2,   b_2) = (0,   0)$, then  $l_2 = \lambda_2$ and  Guichard condition implies that $l_1 = \lambda_2 \cos \varphi(\eta)$ and $l_3 = \lambda_2 \sin \varphi(\eta)$. In order to f\/ind $\varphi$, we use equation~(\ref{E}) with the following indices
\begin{gather*}
h_{13,x_3} + h_{31,x_1} + h_{12} h_{32} = 0.
\end{gather*}
By using the same arguments as in $a)$, we have the following ODE
\begin{gather*}
\varphi_{,\eta \eta} \big[ (N_1)^2 - (N_3)^2 \big] - 2 \varphi_{,\eta} (b_1 N_1 - b_3 N_3) = 0,
\end{gather*}
whose solution will depend on $b_1$ and $b_3$. If $b_1 = b_3$ we have $\varphi$ given by (\ref{eta.dilata21}) and if $b_1 \neq b_3$, the solution is given by~(\ref{eta.dilata22}).

$c)$ The arguments when $(a_3, b_3) = (0,0)$ are the same as in $a)$.
\end{proof}

\begin{remark}
Although our calculation of the symmetry group for the Lam\'e's system has similar techniques to those used by Tenenblat and Winternitz for the intrinsic generalized wave and sine-Gordon equations in \cite{keti3}, we observe that the solutions invariant under the subgroups are quite dif\/ferent. In fact, when we consider the solutions invariant under the translation subgroup in Theorem \ref{teo.nonconstant.trans}, the solutions of
  (\ref{ode.jacobi})  are given by Jacobi elliptic functions that cannot be reduced to elementary functions. Moreover, the only solutions of the Lam\'e's system, which are invariant under the action of the subgroup involving translation and dilations, that depend on all three variables are constant, in contrast to the solutions in \cite{keti3}. The main reason is due to Guichard condition.
\end{remark}

In the next two sections, we will deal with the geometric properties of the Guichard nets and of the conformally f\/lat hypersurfaces associated to the solutions invariant under the 2-dimensional translation subgroup. As we will see in Section~\ref{section5}, these are the solutions that will provide a new class of conformally f\/lat hypersurfaces.

\section{Geometric properties of the Guichard nets}\label{section4}

In this section, we will study the geometric properties of the Guichard nets associated to locally conformally f\/lat hypersurfaces corresponding to the solutions of  the Lam\'e's system $l_i (\xi)$, which are invariant under the translation subgroup. Let $l_1 (\xi)$,   $l_2 (\xi)$,   $l_3 (\xi)$, with $\xi = \sum\limits_{s=1}^3 \alpha_s x_s$ be a~solution of Lam\'e's system. Theorem~\ref{existence} implies that there is a Guichard net $x = (x_1, x_2, x_3) : U \subset \R^3 \rightarrow \R^3$, with a Riemannian metric
\begin{gather}
g = l_1^2 dx_1^2 + l_2^2 dx_2^2 + l_3^2 dx_3^2, \label{netmetric}
\end{gather}
where $U$ is an open set, given by $U = \big\{ (x_1, x_2, x_3) \in \R^3 \, | \, \xi_1 < \xi < \xi_2 \big\}$, where $\xi_1$ and $\xi_2$ are real constants.

\subsection{Level surfaces}\label{section4.1}

In this subsection, we will show that the Guichard nets are foliated by surfaces $\xi = \xi_0$ which are geodesically parallel. Moreover, we will prove that each such surface has f\/lat Gaussian curvature and constant mean curvature that depends on~$\xi_0$.

\begin{definition}\label{geodparalelos}
 Let $M^n$ be a Riemannian manifold and  let $f:M \rightarrow \R$ be a dif\/ferentiable function. The level submanifolds of $f$ are said to be \textit{geodesically parallel} if $|\operatorname{grad} f|$ is a non zero constant,  along each level submanifold.
\end{definition}

We have the following theorem
\begin{theorem}
Let $(U, g)$, $U \subset \R^3$, be a Riemannian manifold with coordinates $(x_1, x_2, x_3)$ and metric $g =  \sum\limits_{s=1}^3 l_s^2 (\xi) dx_i^2$, where $\xi =  \sum\limits_{s=1}^3 \alpha_s x_s$.  Then the level surfaces
\begin{gather*}
P_{\xi_0} = \left\{ (x_1, x_2, x_3) \in U; \ \sum_{s=1}^3 \alpha_s x_s = \xi_0   \right\}, \qquad \text{where} \quad \xi_1 < \xi_0 < \xi_2,
\end{gather*}
endowed with the induced metric, are geodesically parallel. Moreover, each level surface  has flat Gaussian curvature and constant mean curvature (depending on $\xi_0$).
\end{theorem}
\begin{proof} Since at least one  $\alpha_i$ is non zero, we can suppose that  $\alpha_3 \neq 0$ and we parametrize  $P_{\xi_0}$ as
\begin{gather*}
X(x_1, x_2) = \left(x_1, x_2, \frac{\xi_0 - \alpha_1 x_1 - \alpha_2 x_2}{\alpha_3} \right).
\end{gather*}
Then $X_{,x_1} = \left( 1, 0, -\alpha_1 / \alpha_3 \right)$ and $X_{,x_2} = \left( 0, 1, - \alpha_2 / \alpha_3 \right)$. Consequently, the coef\/f\/icients of the induced metric are constant, since $\xi = \xi_0$ in this surface. Therefore the Gaussian curvature is equal to zero.

Consider now the function  $h(x) =  \sum\limits_{i=1}^3 \alpha_i x_i$. Then $P_{\xi_0} = h^{-1}(\xi_0).$ Since  $h$ is constant along~$P_{\xi_0}$, it follows that $\operatorname{grad}h$ is normal to  $P_{\xi_0}$. Moreover,
\begin{gather*}
g(\operatorname{grad}h, \operatorname{grad}h) =  \sum_{j=1}^3 \frac{\alpha_j^2}{l_j^2},
\end{gather*}
which implies that $|\operatorname{grad} h| $ is constant along  $P_{\xi_0}$. It follows from Def\/inition~\ref{geodparalelos} that the level surfaces $h^{-1}(\xi_0)$ are geodesically parallel.

Now we compute the mean curvatures of $P_{\xi_0}$. Given $p \in P_{\xi_0}$,  let  $A: T_p P_{\xi_0} \rightarrow T_p P_{\xi_0}$ be the Weingarten operator, i.e.,  $A \textbf{v} = - \nabla_{\textbf{v}} \left( \frac{\operatorname{grad} h}{|\operatorname{grad} h|} \right)(p)$, where $\nabla$ is the Riemannian connection on  $(U, g)$. Since $|\operatorname{grad} h|$ is constant along  $P_{\xi_0}$, it follows that
\begin{gather*}
A \textbf{v} = - \frac{1}{|\operatorname{grad} h|} \nabla_{\textbf{v}} \operatorname{grad} h (p).
\end{gather*}
Then the mean curvature of $P_{\xi_0}$ is given by
\begin{gather*}
H = - \frac{\Delta h(p)}{|\operatorname{grad} h|}
= \frac{1}{|\operatorname{grad} h|} \sum_{i,k} \frac{\Gamma^k_{ii} (\xi_0) \alpha_k}{l_i^2(\xi_0)} ,
\end{gather*}
where $\Gamma^k_{ij}$ are the Christof\/fel for the connection $\nabla$. Therefore, the mean curvature of  $P_{\xi_0}$ is a~constant depending on~$\xi_0$.
\end{proof}

\subsection{Coordinate surfaces}\label{section4.2}

In this subsection,  we will use the solutions invariant by the group of  translations  to show that the coordinate surfaces of the corresponding Guichard net $(U,g)$ have constant Gaussian curvature. Moreover, the values of these curvatures satisfy an algebraic relation.

\begin{theorem}
Let $(U, g)$, $U \subset \R^3$, be a Riemannian manifold, with coordinates $(x_1, x_2, x_3)$ and metric $g = \sum\limits_{s=1}^3 l_s^2 (\xi) dx_i^2$, with $\xi =  \sum\limits_{s=1}^3 \alpha_s x_s$. Then each coordinate surface of $U \subset \R^3$, $x_i=\operatorname{const}$, endowed
with the induced metric, has constant Gaussian curvature $K_i$. Moreover,
\begin{gather*}
K_1 + K_2 + K_3 = 0.
\end{gather*}
\end{theorem}
\begin{proof}
Since $g$ is given by (\ref{netmetric}), it follows that the metric induced on each coordinate surface, $x_i=\operatorname{const}$, is
\begin{gather*}
g_i = l_j^2 (dx_j)^2 + l_k^2 (dx_k)^2, \qquad i, \, j, \, k \  {\rm distinct},
\end{gather*}
and its Gaussian curvature, $K_i$, is given by
\begin{gather}
K_i = \frac{1}{ l_j l_k} \left( \frac{l_{k,x_i}l_{j,x_i}}{l_i^2} \right). \label{gaussianaparal}
\end{gather}

 Assume that none of the functions $l_i$ is constant and
 $\xi = \alpha_1 x_1 + \alpha_2 x_2 + \alpha_3 x_3$, with $\alpha_i \neq 0$, for all $i$.  In this case, we have $l_{i,\xi} = c_i l_j l_k$, where  $i$, $j$ and $k$ are distinct indices in $\left\{ 1,2,3 \right\}$.
Therefore, it follows from (\ref{gaussianaparal}) that the Gaussian curvature of each coordinate surface is given by $K_i = c_j c_k \alpha_i^2$. Moreover, it follows from  (\ref{alface}) that
\begin{gather*}
 K_1 + K_2 + K_3 = \alpha_1^2 c_2 c_3 + \alpha_2^2 c_3 c_1 + \alpha_1^3 c_1 c_2 = 0.
 \end{gather*}

 If only one of the functions $l_i$  is constant, it follows from Lemma~\ref{ljalphaj}, that, if $l_i$ is constant, then $\alpha_i=0$. Then it follows from~(\ref{gaussianaparal}) that all the curvatures are equal to zero. In fact, $K_i = 0$, since the functions~$l_s$, for all~$s$, do not depend on~$x_i$. Moreover, for $j \neq i$, $K_j = 0$, since $l_i$ is constant. Hence, the sum $\sum\limits_{j=1}^3 K_i = 0$  trivially.
\end{proof}

\section{Conformally f\/lat hypersurfaces}\label{section5}

In this section, we describe the generic conformally f\/lat hypersurfaces associated to the solutions of the Lam\'e's system  invariant under the translation group.

It is known that, any locally generic conformally f\/lat hypersurface, in a 4-dimensional space form, has a metric induced by the Guichard net of the form (see \cite{suyamajer,Suyama2, Suyama3})
\begin{gather}
g = e^{2 P(x)} \big\{ \sin^2 \varphi (x)  dx_1^2 +   dx_2^2 + \cos^2 \varphi (x)  dx_3^2  \big\}, \label{inducedmetrici}
\end{gather}
where $x = (x_1,   x_2,   x_3)$,  or
\begin{gather}
g = e^{2 \tilde{P}(x)} \big\{ \sinh^2 \tilde{\varphi} (x)  dx_1^2 + \cosh^2 \tilde{\varphi} (x)  dx_2^2 +   dx_3^2  \big\}. \label{inducedmetricii}
\end{gather}

Suyama classif\/ied in \cite{Suyama3} the hypersurfaces conformal to the products $M^2 \times I \subset \R^4$ given by Lafontaine in~\cite{Lafontaine}, as the hypersurfaces where $\varphi$ depends only on two variables. Hertrich-Jeromin and Suyama classif\/ied in \cite{suyamajer} the hypersurfaces where $\varphi$ has two  vanishing mixed derivatives. These conformally f\/lat hypersurfaces are associated to the so called \emph{cyclic Guichard nets}, which are characterized by $\varphi_{,x_1 x_2} = \varphi_{,x_2 x_3} = 0$, when $g$ is of the form (\ref{inducedmetrici}) and by $\varphi_{,x_1 x_3} = \varphi_{,x_2 x_3} = 0$, when $g$ is given by~(\ref{inducedmetricii}). Moreover, the authors showed that all the known  cases of conformally f\/lat hypersurfaces, up to now, are associated to cyclic Guichard nets.

We observe that Theorem \ref{uma.constante.dilata} shows that each solution of the Lam\'e's system, which is invariant under the action of the 2-dimensional subgroup involving translations and dilations, depends only on two variables. Therefore, the conformally f\/lat hypersurfaces associated to these solutions are conformal to the products $M^2 \times I \subset \R^4$.

We now consider the conformally f\/lat hypersurfaces associated to the solutions invariant under the translation subgroup. We analyse each case  separately:

$i)$ $\xi = \alpha_1 x_1 + \alpha_2 x_2$.
In this case, we have the solutions $l_1 = \lambda_3 \sinh (\xi + \xi_0)$, $l_2 = \lambda_3 \cosh (\xi + \xi_0)$ and $l_3 = \lambda_3 \neq 0$ (see Theorem \ref{uma.constante}). Then the associated conformally f\/lat hypersurface  has a Guichard net, where the induced metric is given by
\begin{gather}
g = e^{2 P(x_1, x_2, x_3)} \big\{ \sinh^2 (\xi + \xi_0) dx_1^2 + \cosh^2 (\xi + \xi_0)  dx_2^2 + dx_3^2  \big\}. \label{metrichyp}
\end{gather}
The hypersurface is conformal to one of the products considered by Lafontaine in \cite{Lafontaine} that we describe as follows (see \cite{Suyama2, Suyama3} for details): Let $\h^{3}$ be the hyperbolic 3-space, considered as the half space model and  as a subset of $\R^4$, i.e.,
\begin{gather*}
\h^{3}=\big\{\big(y^{1}, y^{2}, y^{3},0\big) : y^{3}>0 \big\} \subset \R^{4}= \big\{\big(y^{1}, y^{2}, y^{3},y^{4}\big) : y^{i} \in \R \big\},
\end{gather*}
with the metric $g_{ij} = \frac{\delta_{ij}}{y_3^2}$. Consider the rotations of the  $y^3$-axis given by
\begin{gather*}
\big(y^{1}, y^{2},y^{3},0\big) \rightarrow \big(y^{1}, y^{2}, y^{3}\cos t, y^{3}\sin t\big),
\end{gather*}
then the hypersurface $M^3 = M^2 \times I$, obtained by the above rotation of a surface of constant curvature $M^2 \subset \h^3$ is a conformally f\/lat hypersurface. One can show that for $g$ given by~(\ref{metrichyp}), the surface $M^2 \subset \h^3$ is a f\/lat surface,  parametrized by lines of curvature whose f\/irst and second fundamental forms are given by
\begin{gather}
I  =  \sinh^2 (\xi + \xi_0) dx_1^2 + \cosh^2 (\xi + \xi_0)  dx_2^2,  \nonumber\\
II  =  \sinh (\xi + \xi_0) \cosh (\xi + \xi_0) \big(dx_1^2 + dx_2^2\big).
 \label{firstffhyp}
\end{gather}

In order to describe the f\/lat surfaces $M^2 \subset \h^3$, we mention a classif\/ication result obtained by the authors in collaboration with Mart\'{\i}nez~\cite{marsanten}. It is well  known that,  on a neighbourhood  of a non-umbilical point, a f\/lat surface
in $\h^3$ can be parametrized by lines of curvature, so that
 the f\/irst and second fundamental forms are given by (for details, see \cite[Theorem~2.4, Corollary~2.7]{keti})
\begin{gather}
I = \sinh^2 \phi(u,v) (du)^2 + \cosh^2 \phi(u,v) (dv)^2, \label{firstff} \\
II = \sinh \phi(u,v) \cosh \phi(u,v) \big((du)^2 + (dv)^2 \big), \label{secondff}
\end{gather}
where $\phi$ is a harmonic function, i.e. $\phi_{uu}+\phi_{vv}=0$. The classif\/ication result is given as follows:

\begin{theorem}[\cite{marsanten}] \label{marsanten}
Let $\Sigma$ be a flat surface in $\h^3$ with a local parametrization, in a neighborhood of a nonsingular and nonumbilic point, such that the first and second fundamental forms are diagonal and given by \eqref{firstff} and \eqref{secondff}, where $\phi$ is a $($Euclidean$)$ harmonic function. Then $\phi$ is linear, i.e., $\phi = au + bv+ c$ if, and only if, $\Sigma$ is locally congruent to either a helicoidal flat surface $($when $(a,b,c) \neq (0, \pm 1, 0))$ or to a ``peach front'' $($when $(a,b,c) = (0, \pm 1, 0))$.
\end{theorem}

Helicoidal surfaces arise as a natural generalization of rotational surfaces. They are invariant by a helicoidal group of isometries, i.e., given an axis, we consider a translation along this axis composed with a rotation around it. In the half space model of~$\h^3$, up to isometries, we can consider the $y_3$-axis, which enables us to write the helicoidal group, relative to this axis, as the composition
\begin{gather*}
  h_t = \left(
  \begin{matrix}
    e^{  \beta t} & 0 & 0 \\
    0 & e^{ \beta t} & 0 \\
    0 & 0 & e^{\beta t}
  \end{matrix}
  \right)
  \left(
  \begin{matrix}
    \cos \alpha t & - \sin \alpha t & 0 \\
    \sin \alpha t & \cos \alpha t & 0 \\
    0 & 0 & 1
  \end{matrix}
  \right),
\end{gather*}
of a rotation around the $y_3$-axis with \textit{angular pitch} $\alpha$ with a hyperbolic translation of \textit{ratio}~$\beta.$ The ``peach front'' is a special case of f\/lat surfaces that is not helicoidal. Details about this surface can be found in \cite{flatfront1}.

The study of f\/lat surfaces in hyperbolic 3-space has received a lot of attention in the last few years, mainly because Galv\'ez, Mart\'{\i}nez and Mil\'an have shown in~\cite{galvez} that f\/lat surfaces in the hyperbolic 3-space admit a Weierstrass representation formula in terms of meromorphic data as in the theory  of minimal surfaces in $\R^3$. Namely, if $\psi : M^2 \rightarrow \h^3$ is a surface in~$\h^3$, for any $p \in M^2$, there exist $G(p),   G^{*}(p) \in \mathbb{C}_{\infty}$ distinct points in the ideal boundary of $\h^3$ such that the oriented normal geodesic at $\psi(p)$ is the geodesic in $\mathbb{H}^3$ starting from $G^{*}(p)$ towards~$G(p)$. The maps $G, G^{*}: \Sigma \rightarrow \mathbb{C}_{\infty}$ are called the \textit{hyperbolic Gauss maps} of~$\psi$. It is proved in~\cite{galvez} that, for f\/lat surfaces, they are holomorphic when one considers~$\mathbb{C}_{\infty}$ as the Riemann sphere and~$M^2$ has a~complex structure induced by the second fundamental form. Conversely, given two holomorphic functions~$G$ and~$G^{*}$, with $G \neq G^{*}$, one can recover a f\/lat immersion of a surface in~$\h^3$ (for more details see also~\cite{CorMilMart, flatfront1, flatfront2}). This representation formula was the main tool to obtain Theorem~\ref{marsanten}.

With the previous results we conclude that

\begin{theorem} \label{characterization.helicoidal}
Let $l_i(\xi)$ be solutions of the Lam\'e's system, where $\xi = \alpha_1 x_1 + \alpha_2 x_2$. Then the associated conformally flat hypersurfaces are conformal to the product, $M^2 \times I$,  where $M^2$ is locally congruent to either a helicoidal flat surface in $\h^3$ or the ``peach front''.
\end{theorem}

\begin{proof}
When $\xi = \alpha_1 x_1 + \alpha_2 x_2$, it follows from Theorem~\ref{uma.constante} that  the solution of Lam\'e's system  is $l_1 = \lambda_3 \sinh (\xi + \xi_0)$, $l_2 = \lambda_3 \cosh (\xi + \xi_0)$ and $l_3 = \lambda_3 \neq 0$ and the corresponding conformally f\/lat hypersurface $M^3$ has a metric $g$ given by~(\ref{metrichyp}). Then $M^3$ is conformal to the product $M^2 \times I$, where $M^2$ is a f\/lat surface in $\h^3$ with  fundamental forms given by~(\ref{firstffhyp}). It follows from  Theorem~\ref{marsanten} that $M^2$ is locally congruent to either a helicoidal f\/lat surface in $\h^3$ or to the ``peach front''. \end{proof}

$ii)$ $\xi = \alpha_1 x_1 + \alpha_3 x_3$.
In this case, we have the solution of Lam\'e's system , $l_1 = \lambda_2 \sin (\xi + \xi_0)$, $l_2 = \lambda_2\neq 0 $ and $l_3 = \lambda_2 \cos (\xi + \xi_0)$ (see Theorem \ref{uma.constante} b)). Then the associated conformally f\/lat hypersurface $M^3$ has a Guichard net, whose induced metric is given by
\begin{gather*}
g = e^{2 P(x_1, x_2, x_3)} \big\{ \sin^2 (\xi + \xi_0) dx_1^2 +   dx_2^2 + \cos^2 (\xi + \xi_0) dx_3^2  \big\}. 
\end{gather*}
The hypersurface $M^3$ is conformal to another class of products $M^2 \times I$ (see \cite{Suyama2,Suyama3}). Namely, let $\mathbb{S}^{3} \subset \R^4$ be the canonical 3-sphere, then
$M^2 \times I = \left\{ tp:0<t<\infty, p\in M^2 \subset \mathbb{S}^3 \right\}$, is a~conformally f\/lat hypersurface, where $M^2$ is a~surface with constant curvature in $S^3$. In our case, $M^3$ is conformal to the  product $M^2 \times I$, where the surface $M^2 \subset \mathbb{S}^3$ is a f\/lat surface,  parametrized by lines of curvature, whose f\/irst and second fundamental forms are  given by
\begin{gather}
I =  \sin^2 (\xi + \xi_0) dx_1^2 + \cos^2 (\xi + \xi_0)  dx_3^2,  \nonumber\\
II  =  \sin (\xi + \xi_0) \cos (\xi + \xi_0) \big(dx_1^2 - dx_3^2\big).
 \label{firstffsphere}
\end{gather}

The geometry of these surfaces in $\mathbb{S}^3$ is being studied and it will appear in another paper.

$iii)$ $\xi = \alpha_1 x_1 + \alpha_2 x_2 + \alpha_3 x_3$, $\alpha_i \neq 0$, for all $i$.
In this case, we will show that the solutions~$l_i (\xi)$ of the Lam\'e's system give rise to a new class of conformally f\/lat hypersurfaces, according to the following theorem:

\begin{theorem} \label{newclass}
Let $M^3$ be a conformally flat hypersurface in a space form $M^4_K$, associated to a~solution of Lam\'e's system $l_i (x_1, x_2, x_3) = l_i (\xi)$, with $\xi = \sum\limits_{s=1}^3 \alpha_s x_s$ and $\alpha_s \neq 0$, for all $s$, given in terms of elliptic functions by \eqref{ode.jacobi}--\eqref{l3.jacobi}.  Then its first fundamental form  $g$ is given by
\begin{gather}
g = e^{2P(x)} \big\{ \cos^2 \varphi (\xi) (dx_1)^2 + (dx_2)^2 + \sin^2 \varphi (\xi) (dx_3)^2 \big\},	
	\label{metricacosgeral}
\end{gather}
where $\varphi$ satisfies,
\begin{gather}
\varphi_{,\xi}^2 = c (a \cos^2 \varphi - b ) , \label{fiduplacos}
\end{gather}
or $g$ is given by
\begin{gather}
g = e^{2\tilde{P}(x)} \big\{ \sinh^2 \tilde{\varphi} (\xi) (dx_1)^2 + \cosh^2 \tilde{\varphi} (\xi) (dx_2)^2 + (dx_3)^2 \big\}
	\label{metricacoshgeral}
\end{gather}
where $\tilde{\varphi}$ satisfies
\begin{gather}
\tilde{\varphi}_{,\xi}^2 = c \big(b \cosh^2 \tilde{\varphi} - b \big) . \label{fiduplacosh}
\end{gather}
where $a, b,c \in \R\setminus \{0\}$, $P(x)$ and $\tilde{P}(x)$ are differentiable functions that depend on $l_s$ and  $M^4_K$. In both cases, $\xi \in I \subset \R$, where $I$ is an open interval such that $g$ is positive definite.
\end{theorem}
\begin{proof}
Guichard condition (\ref{relacaoguichard}) implies that we may consider
\begin{gather}
  l_1 = l_2 \cos \varphi, \label{l1cos} \\
  l_3 = l_2 \sin \varphi. \label{l3cos}
 \end{gather}
It follows from Theorem \ref{existence} that the metric is given by (\ref{metricacosgeral}). In order to obtain the expression for the derivative of $\varphi$ with respect to $\xi$, we consider
\begin{gather*}
 l_{1,\xi} = l_{2,\xi} \cos \varphi - l_2 \varphi_{,\xi} \sin \varphi.
\end{gather*}
Since $\alpha_s \neq 0$ for all $s$, the functions $l_i$ are given as in Theorem~\ref{teo.nonconstant.trans}, by (\ref{ode.jacobi})--(\ref{l3.jacobi}). Hence, using~(\ref{derivadalame}), we have that
\begin{gather*}
c_1 l_2 l_3 = c_2 l_1 l_3 \cos \varphi - \varphi_{,\xi} l_3,
\end{gather*}
for $c_1,\,c_2\in\R\setminus \{0\}$, which implies
\begin{gather}
\varphi_{,\xi} = l_2 \big(c_2 \cos^2 \varphi - c_1 \big). \label{1derivadafi}
\end{gather}
By taking the derivative again, it follows from (\ref{l1cos})--(\ref{1derivadafi}) and (\ref{derivadalame}) that
\begin{gather*}
\varphi_{,\xi \xi}  =  l_{2, \xi} \big(c_2 \cos^2 \varphi - 1 \big) - 2 c_2 l_2 \cos \varphi \sin \varphi \varphi_{,\xi}
 =  \frac{1}{l_2} \big[ l_{2,\xi} \varphi_{,\xi} - 2 c_2 (l_2 \cos \varphi) (l_2 \sin \varphi) \varphi_{,\xi}  \big] \\
\hphantom{\varphi_{,\xi \xi}}{}
 =  \frac{1}{l_2} \big[  l_{2,\xi} \varphi_{,\xi} - 2c_2 l_1 l_3 \varphi_{,\xi} \big]
 =  \frac{1}{l_2} \big[  l_{2,\xi} \varphi_{,\xi} - 2 l_{2,\xi} \varphi_{,\xi} \big]
 =  - \frac{l_{2,\xi} \varphi_{,\xi}}{l_2}.
\end{gather*}
Therefore,
\begin{gather}
 \varphi_{,\xi} l_2 = c,
\label{fil2}
\end{gather}
where $c\in\R\setminus \{0\}$, since  $\varphi_{,\xi} \neq 0$. Then, multiplying~(\ref{1derivadafi}) by  $\varphi_{,\xi}$ and using (\ref{fil2}) we have that $\varphi_{,\xi}^2 = c \big(c_2 \cos^2 \varphi - c_1 \big)$, i.e., (\ref{fiduplacos}) holds. The proof of the second part of the theorem is analogous, when we consider $l_2 = l_3 \cosh \varphi$ and $l_1 = l_3 \sinh \varphi$.
\end{proof}

\begin{corollary} \label{cor.newclass}
Let $M^3 \subset M^4$ be a conformally flat hypersurface associated to the solutions of Lam\'e's system $l_i(\xi)$ with $\xi =  \sum\limits_{s=1}^3 \alpha_s x_s$ and $\alpha_s \neq 0$ for all $s$, given in terms of elliptic functions by \eqref{ode.jacobi}--\eqref{l3.jacobi}. Then the Guichard net of $M^3$ is not cyclic.
\end{corollary}

\begin{proof}
It follows from Theorem~\ref{newclass} that the f\/irst fundamental form of $M^3$ is given by~(\ref{metricacosgeral}), where $\varphi(\xi)$ satisf\/ies (\ref{fiduplacos}) or by (\ref{metricacoshgeral}) where $\tilde{\varphi} (\xi)$ satisf\/ies~(\ref{fiduplacosh}), $\xi \in I \subset \R$. In the f\/irst case, (\ref{fiduplacos}) implies that $\varphi_{,\xi\xi} = \lambda \sin 2 \varphi$, $\lambda \neq 0$ and in the second case, (\ref{fiduplacosh}) implies that $\tilde{\varphi}_{,\xi\xi} = \lambda \sin 2 \tilde{\varphi}$, $\lambda \neq 0$. In either case, $\varphi_{,x_i x_j} = \alpha_i \alpha_j \varphi_{,\xi\xi} \neq 0$, $i \neq j$ and $\xi \in I$. Therefore, the Guichard net of~$M^3$ is not cyclic.
\end{proof}

We observe that, as a consequence of  the results of Hertrich-Jeromin and Suyama, the surfaces~$M^3$ of Corollary~\ref{cor.newclass} provide a new class of conformally f\/lat hypersurfaces.

It is important to observe that Hertrich-Jeromin and Suyama in \cite{suyamajer2} have  independently considered Guichard nets with the ansatz on the function $\varphi$ such that $\varphi (x_1, x_2, x_3) = \varphi ( a x_1 + b x_2 + c x_3)$. They investigated the geometric properties of these Guichard nets, that they called  Bianchi-type Guichard nets, as well as the new class of associated conformally f\/lat hypersurfaces.

\appendix

\section{Appendix}\label{appendix}

\begin{proof}[Proof of Theorem \ref{teo.generator}.]
The inf\/initesimal generator associated to the symmetry group is written as
in (\ref{infgeneratorlame}). The functions  $\xi^i$, $\eta^i$, $\phi^{ij}$ will be obtained by solving the determining equations that arise when we apply the f\/irst prolongation formula
 \begin{gather*}
{\rm pr}^{(1)} V = V +  \sum_{i,k} D_k\big(\eta^i\big) \frac{\partial}{\partial l_{i,x_k}} +  \sum_{i,j,k \, i \neq j} D_k\big(\phi^{ij}\big)  \frac{\partial}{\partial h_{ij,x_k}}\\
\hphantom{{\rm pr}^{(1)} V =}{}
-  \sum_{i,k,r} D_k(\xi^r) l_{i,x_r} \frac{\partial}{\partial l_{i,x_k}} -  \sum_{i,j,k,r} D_k(\xi^r) h_{ij,x_r} \frac{\partial}{\partial h_{ij,x_k}},
\end{gather*}
with
\begin{gather*}
D_i = \frac{\partial}{\partial x_i} +   \sum_j l_{j,x_i} \frac{\partial}{\partial l_i} +   \sum_{j,l} h_{jl,x_i} \frac{\partial}{\partial h_{jl}},
\end{gather*}
on each equation of the system, i.e., when we consider the inf\/initesimal criterion~(\ref{criterion}). In order to avoid any functional dependence,  the following substitutions will be considered
\begin{gather}	
	l_{i,x_j} = h_{ij}l_j   ,  \qquad i  \neq j, \label{subst.lij} \\
	l_{i,x_i} = - \ei \ej h_{ji} l_j - \ei \ek h_{ki} l_k \label{subst.lii}  \\
	h_{ij,x_k} = h_{ik} h_{kj},  \label{subst.hijk} \\
	h_{ij,x_j} = - h_{ji,x_i} - h_{ik} h_{jk}  , \qquad  i < j, \label{subst.hijj} \\
	h_{ij,x_i} = - \ei \ej h_{ji,x_j} - \ei \ek h_{ki} h_{kj}  , \qquad i < j. \label{subst.hiji}
\end{gather}

Fixing $i$, $j$ and $k$, distinct indices, we start applying  ${\rm pr}^{(1)} V$ to  equation (\ref{D}). Then the inf\/initesimal criterion (\ref{criterion}), gives $\phi_{(k)}^{ij} - \phi^{ik}h_{kj} - h_{ik} \phi^{kj} = 0$, using the prolongation formula, we get
\begin{gather}
\phi_{,x_k}^{ij} + \sum_r \phi^{ij}_{,l_r} l_{r,x_k} + \sum_{r,s} \phi ^{ij}_{,h_{rs}} h_{rs,x_k}
 - \sum_t \left( \xi_{,x_k}^t + \sum_r \xi^t_{,l_r} l_{r,x_k} + \sum_{r,s} \xi^t_{,h_{rs}} h_{rs,x_k} \right) h_{ij,x_t}  \nonumber\\
 \qquad{} - \phi^{ik} h_{kj} - h_{ik} \phi^{kj} = 0 .
  \label{1prolong}
\end{gather}
For $i<j$, we apply the substitutions (\ref{subst.lij})--(\ref{subst.hiji}) and we analyse each term of (\ref{1prolong}) as follows
\begin{gather}
 \sum_r \phi^{ij}_{,l_r} l_{r,x_k} =  \sum_{r \neq k} \phi^{ij}_{,l_r} h_{rk} l_k - \phi^{ij}_{,l_k} \left( \ek \ej h_{jk} l_j + \ek \ei h_{ik} l_i  \right), \label{1prolong.1soma}
\\
 \sum_{r,s} \phi^{ij}_{,h_{rs}} h_{rs,x_k}  = \! \sum_{r \neq k   , \, s \neq k} \!\phi^{ij}_{,h_{rs}} h_{rk} h_{ks}
 + \sum_{s < k}\! \phi^{ij}_{,h_{ks}} h_{ks,x_k}
 - \sum_{s > k} \!\phi^{ij}_{,h_{ks}} \big(  \ek \varepsilon_s h_{sk,x_s} + \ek \varepsilon_m h_{mk}h_{ms} \big) \nonumber\\
\hphantom{\sum_{r,s} \phi^{ij}_{,h_{rs}} h_{rs,x_k}  =}{}
 - \!\sum_{r < k} \phi^{ij}_{,h_{rk}} \left(  h_{kr,x_r} + h_{rn} h_{kn} \right) + \sum_{r > k} \phi_{,h_{rk}}^{ij} h_{rk,x_k},\label{1prolong.2soma}
\\
 \sum_t \left( \xi^t_{,x_k} + \sum_r \xi^t_{,l_r}l_{r,x_k} + \sum_{r,s} \xi^t_{,h_{rs}} h_{rs,x_k} \right) h_{ij,x_t} \nonumber\\
 \qquad{}
 =  C^k_k h_{ik}h_{kj} - C^j_k \big( h_{ji,x_i} + h_{ik} h_{jk} \big)
  - C^i_k \big( \ei \ej h_{ji,x_j} + \ei \ek h_{ki} h_{kj} \big) ,\label{1prolong.3soma}
\end{gather}
where the coef\/f\/icients $C^t_k$ are given by
  \begin{gather}
	C^t_k  =   \xi^t_{,x_k} + \sum_{r \neq k} \xi^t_{,l_r} h_{rk} l_k - \xi^t_{,l_k} \big( \ek \ej h_{jk} l_j + \ek \ei h_{ik} l_i \big) + \sum_{r \neq k, \, s \neq k} \xi^t_{,h_{rs}}h_{rk}h_{ks} \nonumber\\
\hphantom{C^t_k  =}{}
+ \sum_{s < k} \xi^t_{,h_{ks}} h_{ks,x_k}	  -  \sum_{s>k}\xi^t_{,h_{ks}} \big( \ek \varepsilon_s h_{sk,x_s} + \ek \varepsilon_m  h_{mk} h_{ms} \big)
\nonumber\\
\hphantom{C^t_k  =}{}
	   - \sum_{r<k} \xi^t_{,h_{rk}} \big( h_{kr,x_r} + h_{rn} h_{kn} \big) + \sum_{r>k} \xi^t_{,h_{rk}} h_{rk,x_k},
  \label{ctkfunction}
\end{gather}
 and the indices $m$ and $n$ are such that $ \{ k,  s,  m  \}$, $ s>k$ and $ \{ k,  r,  n  \}$, $r<k$ are two sets of three distinct numbers.

Now we analyse the coef\/f\/icients of equation (\ref{1prolong}), considering  (\ref{1prolong.1soma})--(\ref{1prolong.3soma}). By equating to zero the coef\/f\/icients of the products $h_{ji,x_j} h_{ks,x_k}$, with  $k>s$, we obtain  $\xi^i_{,h_{ks}} =0 $. Analogously, for the coef\/f\/icients of  $h_{ji,x_j} h_{sk,x_s}$, with $k<s$, we obtain $\xi^i_{,h_{ks}} = 0$. 	
This implies that
\begin{gather*}
\xi^i_{,h_{ks}} = 0, \quad \forall \, s, \ \  s \neq k, \qquad {\rm i.e.} \quad \xi^i_{,h_{kj}}= \xi^i_{,h_{ki}}=0.
\end{gather*}
Similarly, from the coef\/f\/icients of  $h_{ji,x_j} h_{kr,x_k}$, $r<k$ and $h_{ji,x_j} h_{rk,x_k}$, with $r>k$, we obtain
\begin{gather*}
\xi^i_{,h_{rk}} = 0, \quad \forall \, r, \ \ r \neq k, \qquad {\rm i.e.} \quad \xi^i_{,h_{jk}}= \xi^i_{,h_{ik}}=0,
\end{gather*}
where $i,   j,   k \in  \{ 1,   2,   3  \}$ are distinct and $i<j$. By analysing  the coef\/f\/icients of  $h_{ji,x_i}h_{ks,x_k}$ with $k>s$ and $h_{ji,x_i}h_{sk,x_s}$ with $k<s$, we obtain
\begin{gather*}
\xi^j_{,h_{ks}} = 0, \quad \forall \, s, \ \  s \neq k, \qquad {\rm  i.e.} \quad \xi^j_{,h_{ki}}= \xi^j_{,h_{kj}}=0.
\end{gather*}
Similarly, the coef\/f\/icients of  $h_{ji,x_i}h_{kr,x_r}$, with $k>r$, and $h_{ji,x_i}h_{rk,x_k}$, with $k<r$,  lead to
\begin{gather*}
\xi^j_{,h_{rk}} = 0, \quad \forall \, r, \ \  r \neq k, \qquad {\rm i.e.} \quad \xi^j_{,h_{ik}}= \xi^j_{,h_{jk}}=0.
\end{gather*}
Since $i,   j,   k \in  \{ 1,   2,   3  \}$ are distinct and arbitrary indices, with $i<j$, we conclude that  $\xi^m_{,h_{st}}= 0$ for any indices $m$, $s$ and $t$, $s \neq t$, i.e., $\xi^m$ depends only on  $x$ and $l$.  Therefore,  the expression of~$C^t_k$ given in~(\ref{ctkfunction}) reduces to
\begin{gather*}
C^t_k = \xi^t_{,x_k} +  \sum_{r \neq k} \xi^t_{,l_r} h_{rk} l_k - \xi^t_{,l_k}  ( \ek \ej h_{jk} l_j + \ek \ei h_{ik} l_i  ),
\end{gather*}
that can be rewritten as
\begin{gather}
	C^t_k= \xi^t_{,x_k} +  \sum_{r \neq k} \big( \xi^t_{,l_r} l_k - \xi^t_{,l_k} \varepsilon_r \ek l_r \big)h_{rk}.
	\label{novoct}
\end{gather}

From  (\ref{1prolong.2soma}), we have that the coef\/f\/icients of  $h_{ks,x_k}$, with $s<k$, and the coef\/f\/icients of $h_{sk,x_s}$, with $s>k$, lead to  $\phi^{ij}_{,h_{ks}} = 0$, $\forall\, s \neq k$,  i.e.\ $\phi^{ij}_{,h_{ki}} = \phi^{ij}_{,h_{kj}} = 0$.

Considering (\ref{1prolong.3soma}), the coef\/f\/icients of $h_{ji,x_i}$ and $h_{ji,x_j}$ imply that $C^j_k=0$ and $C^i_k=0$, respectively. Since $i<j$ and $i, j,  k \in \{ 1,  2,  3 \}$ are arbitrary and distinct, we conclude that $C^i_k = C^j_k = 0$, for all $i$, $j$, $k$ distinct indices.

Since $\xi^m$ does not depend on $h_{st}$, the analysis of  (\ref{novoct}) gives us the following system
\begin{gather*}
	\xi^m_{,x_k} = 0, \qquad  \forall\, m \neq k, \\
	\xi^m_{,l_r} l_k - \varepsilon_r \ek \xi^m_{,l_k} l_r  =  0, \qquad \forall \, r \neq k.
\end{gather*}
The f\/irst equation of this system says  that $\xi^m$ depends only on $x_m$ and $l$. By solving the characteristic system for the second equation, we have that $\xi^m$ depends on $x_m$ and a variable $\zeta =\ei l_i^2 + \ej l_j^2 + \ek l_k^2$. However,  Guichard condition implies that $\zeta \equiv 0$, hence $\xi^m$ does not depend on $l_s$, for all $s$.

Summarising the conclusions of this f\/irst part of the proof,  we have that
\begin{gather*}
\phi^{st} = \phi^{st}  ( h_{st}, h_{ts}, x, l  ), \qquad {\rm and} \qquad \xi^m = \xi^m (x_m). 
\end{gather*}

We now consider equation (\ref{B}). By applying the prolongation ${\rm pr}^{(1)}V$ to (\ref{B}), we have that $\eta^i_{(j)} - \phi^{ij}l_j - h_{ij} \eta^j =0$, which implies,
\begin{gather}
	\eta^i_{,x_j} +  \sum_r \eta^i_{,l_r} l_{r,x_j} + \sum_{r,s} \eta^i_{,h_{rs}} h_{rs,x_j} - \xi^j_{,x_j} l_{i,x_j} - \phi^{ij} l_j - h_{ij} \eta^j = 0.
	\label{2prolong}
\end{gather}
Observe that by applying the substitution  (\ref{subst.lij}), we have
\begin{gather*}
 \sum_r \eta^i_{,l_r} l_{r,x_j} = \sum_{r \neq j} \eta^i_{,l_r} h_{rj} l_j - \eta^i_{,l_j}  ( \ej \ei h_{ij} l_i + \ej \ek h_{kj} l_k  ).
\end{gather*}
Moreover, by applying the substitutions (\ref{subst.hijk}), (\ref{subst.hijj}) and (\ref{subst.hiji}) we have
\begin{gather*}
\sum_{r, \, s} \eta^i_{,h_{rs}} h_{rs,x_j} = \sum_{r \neq j, \, s \neq j} \eta^i_{,h_{rs}} h_{rj} h_{js} + \sum_{s<j} \eta^i_{,h_{js}} h_{js,x_j} - \sum_{s>j} ( \ej \varepsilon_s h_{sj,x_s} + \ej \varepsilon_m h_{mj} h_{ms}  ) \\
\hphantom{\sum_{r, \, s} \eta^i_{,h_{rs}} h_{rs,x_j} =}{}
  -  \sum_{r<j} \eta^i_{,h_{rj}}  ( h_{jr,x_r} + h_{rn} h_{jn}  ) + \sum_{r>j} \eta^i_{,h_{rj}} h_{rj,x_j}.
\end{gather*}
Therefore, by considering in (\ref{2prolong}), the coef\/f\/icients of $h_{js,x_j}$, with  $s<j$, and $h_{sj,x_s}$, with $s>j$, we conclude that  $\eta^i_{,h_{js}} = 0$. Similarly, the analysis of the coef\/f\/icients of~$h_{jr,x_r}$, with~$r<j$, and~$h_{rj,x_j}$, with~$r>j$, imply that   $\eta^i_{,h_{rj}} = 0$.	Hence, we conclude that
\begin{gather*}
\eta^i_{,h_{jt}} = \eta^i_{,h_{tj}} = 0, \qquad \forall \, t \neq j.	
\end{gather*}
Since $i$ and $t \neq j$ are arbitrary, we conclude that $\eta^m$ does not depend on  $h_{st}$, for any indices,~$m$,~$s$ and~$t$ with $s \neq t$. Consequently, (\ref{2prolong}) reduces to
\begin{gather}
	\eta^i_{,x_j} + \big( \eta^i_{,l_i} l_j - \ei \ej \eta^i_{,l_j}l_i - \xi^j_{,x_j} l_j - \eta^j \big ) h_{ij}
	+ \big( \eta^i_{,l_k} l_j - \ej \ek \eta^i_{,l_j} l_k \big) h_{kj} - \phi^{ij} l_j = 0. \label{2prolong.reduzida}
\end{gather}
Since $\phi^{ij}$ depends only on  $x$, $l$, $h_{ij}$ and $h_{ji}$, we obtain from (\ref{2prolong.reduzida}) the following system
\begin{gather}
\eta^i_{,l_k} l_j - \ej \ek \eta^i_{,l_j} l_k = 0, \label{prolong.resto.1} \\
\eta^i_{,x_j} + \big( \eta^i_{,l_i} l_j - \ei \ej \eta^i_{,l_j}l_i - \xi^j_{,x_j} l_j - \eta^j \big) h_{ij} - \phi^{ij} l_j = 0. \label{2prolong.resto.2}	
\end{gather}
By solving the characteristic system for (\ref{prolong.resto.1}), we have  that $n^i=n^i(x,l_i)$. By taking derivatives of (\ref{2prolong.resto.2}) with respect to~$h_{ji}$ we get
\begin{gather}
\phi^{ij}_{,h_{ji}} = 0.	
	\label{conclusao.fi.ji}
\end{gather}
On the other hand, by taking the derivatives of  (\ref{2prolong.resto.2}) twice with respect to $h_{ij}$, we obtain
\begin{gather}
	\phi^{ij}_{,h_{ij}h_{ij}} = 0. \label{conclusao.fi.ij}
\end{gather}
Consequently, it follows from (\ref{conclusao.fi.ji}) and  (\ref{conclusao.fi.ij}) that $\phi^{ij}$ is given by
\begin{gather}
	\phi^{ij} = A^{ij}(x,l)h_{ij} + B^{ij} (x,l).
	\label{conclusao.fi}
\end{gather}
Therefore, (\ref{1prolong}) reduces to
\begin{gather*}
	\phi^{ij}_{,x_k} +  \sum_r \phi^{ij}_{,l_r} l_{r,x_k} + A^{ij} h_{ij,x_k} - \xi^k_{,x_k} h_{ij,x_k} - \big( A^{ik} h_{ik} + B^{ik} \big) h_{kj} - \big( A^{kj}h_{kj} + B^{kj} \big) h_{ik} = 0 .
\end{gather*}
By considering the substitutions (\ref{subst.lij})--(\ref{subst.hiji}), this equation reduces to
\begin{gather*}
	\phi^{ij}_{,x_k} + \phi^{ij}_{,l_i} h_{ik} l_k + \phi^{ij}_{,l_j} h_{jk} l_k -  \phi^{ij}_{,l_k}  ( \ek \ej h_{jk}l_j + \ek \ei h_{ik}l_i  )   \\
\qquad{} 	+ A^{ij} h_{ik} h_{kj} - \xi^k_{,x_k} h_{ik} h_{kj} - A^{ik} h_{ik} h_{kj} - B^{ik}h_{kj} - A^{kj}h_{ik}h_{kj} - B^{kj}h_{ik}  = 0,
	\end{gather*}
which can be rewritten as
\begin{gather*}
	B^{ij}_{,x_k} + A^{ij}_{,x_k}h_{ij} + \big( B^{ij}_{,l_i}l_k - \ek \ei B^{ij}_{,l_k}l_i - B^{kj}  \big) h_{ik} + \big( B^{ij}_{,l_j} l_k - \ek \ej B^{ij}_{,l_k} l_j  \big)h_{jk} \\
\qquad {}- B^{ik}h_{kj} +	 \big( A^{ij}_{,l_i}l_k - \ek \ei A^{ij}_{,l_k} l_i \big) h_{ik}h_{ij}
	 + \big( A^{ij}_{,l_j} l_k - \ek \ej A^{ij}_{,l_k} l_j \big)h_{ij}h_{jk}\\
\qquad {}
	+ \big( A^{ij} - \xi^k_{,x_k} - A^{ik} - A^{kj} \big)h_{ik}h_{kj} =0.
\end{gather*}
It follows from the coef\/f\/icients of $h_{kj}$ that $B^{ik}=0$. The permutation of the indices $i$, $j$ and $k$ leads to
\begin{gather}
B^{st}=0, \qquad \forall \, s,   t, \ \  s \neq t.  \label{conclusao.B}
\end{gather}
By equating to zero the coef\/f\/icients of $h_{ik}h_{kj}$ and $h_{ij} h_{jk}$, the following system is obtained
\begin{gather*}
A^{ij}_{,l_i}l_k - \ek \ei A^{ij}_{,l_k} l_i  = 0, \qquad
A^{ij}_{,l_j} l_k - \ek \ej A^{ij}_{,l_k}	l_j  = 0,
\end{gather*}
where we solve the characteristic system to conclude that $A^{ij}$ depends only on $x$. On the other hand, the coef\/f\/icient of $h_{ij}$ implies that $A^{ij}$ does not depend on  $x_k$, therefore  $A^{ij} = A^{ij}  (x_i, x_j  )$. Considering the coef\/f\/icient of  $h_{ik}h_{kj}$, we obtain the following equation
\begin{gather}
	A^{ij} - \xi^k_{,x_k} - A^{ik} - A^{kj} = 0.
	\label{conclusao.A2}
\end{gather}
Therefore, equation (\ref{2prolong.reduzida}) reduces to
\begin{gather*}
	\eta^i_{,x_j} + \big( \eta^i_{,l_i}l_j - \xi^j_{,x_j} l_j - \eta^j - A^{ij} l_j \big)h_{ij} = 0.
\end{gather*}
Since $\eta^i$ does not depend on $h_{ij}$, we must have
\begin{gather}
\eta^i_{,x_j} = 0, \nonumber\\ 
\eta^i_{,l_i}l_j - \xi^j_{,x_j} l_j - \eta^j - A^{ij} l_j = 0, \label{eta.l}
\end{gather}

By applying ${\rm pr}^{(1)}V$ to equation  (\ref{E}), we have $\phi^{ij}_{(j)} + \phi^{ji}_{(i)} + \phi^{ik}h_{kj} + h_{ik}\phi^{kj} = 0$, which implies that
\begin{gather*}
	\phi^{ij}_{,x_j} + \phi^{ij}_{,h_{ij}} h_{ij,x_j} - \xi^j_{,x_j}h_{ij,x_j} + \phi^{ji}_{,x_i} + \phi^{ji}_{,h_{ji}} h_{ji,x_i} - \xi^i_{,x_i}h_{ji,x_i} + + \phi^{ik}h_{kj} + h_{ik}\phi^{kj} = 0.
\end{gather*}
Considering the substitution (\ref{subst.hijj}), for $i<j$, we obtain
\begin{gather*}
	\phi^{ij}_{,x_j} + \phi^{ji}_{,x_i} + \big( \xi^j_{,x_j} - A^{ij} + A^{ji} - \xi^i_{,x_i} \big) h_{ji,x_i} + \big( \xi^j_{,x_j} - A^{ij} + A^{ik} + A^{kj} \big)h_{ik}h_{jk}  = 0.
\end{gather*}
Then, the coef\/f\/icient of  $h_{ji,x_i}$ leads to
\begin{gather}
	\xi^j_{,x_j} - A^{ij} + A^{ji} - \xi^i_{,x_i} = 0.
	\label{conclusao.Aexi1}
\end{gather}

By applying the prolongation ${\rm pr}^{(1)}V$ to (\ref{F}) we have $\ei \phi^{ij}_{(i)} + \ej \phi^{ji}_{(j)} + \ek \phi^{ki}h_{kj} + \ek \phi^{kj}h_{ki} = 0$,
which implies
\begin{gather*}
	\ei \big( \phi^{ij}_{,x_i} + \phi^{ij}_{,h_{ij}} h_{ij,x_i} - \xi^i_{,x_i} h_{ij,x_i} \big) + \ej \big( \phi^{ji}_{,x_j} + \phi^{ji}_{,h_{ji}} h_{ji,x_j} - \xi^j_{,x_j} h_{ji,x_j} \big) \\
\qquad{} + \ek \big( \phi^{ki} h_{kj} + \phi^{kj}h_{ki} \big) = 0.
\end{gather*}
  Considering the substitution (\ref{subst.hiji}) with  $i<j$, we obtain
\begin{gather*}
	\ei \phi^{ij}_{,x_j} + \ej \phi^{ji}_{,x_j} + \ej \big(   \xi^i_{,x_i} - A^{ij} + A^{ji} - \xi^j_{,x_j} \big) h_{ji,x_j} + \ek \big(   \xi^i_{,x_i} - A^{ij} + A^{ki} + A^{kj} \big) h_{ki} h_{kj} = 0.
\end{gather*}
From the coef\/f\/icient of $h_{ji,x_j}$, we get
\begin{gather}
	A^{ij} - \xi^i_{,x_i} - A^{ji} + \xi^j_{,x_j}= 0.
	\label{conclusao.Aexi2}
\end{gather}
Therefore, it follows from (\ref{conclusao.Aexi1}) and  (\ref{conclusao.Aexi2}) that
\begin{gather}
A^{ji}=A^{ij}.	
	\label{a.simetrico}
\end{gather}
Consequently, both equations imply that $\xi^i_{,x_i} = \xi^j_{,x_j}$, which enables us to conclude that
\begin{gather}
	\xi^m = a x_m + a_m, \qquad \forall \, 1 \leq m \leq 3,
	\label{conclusao.xi.final}
\end{gather}
where $a$ and $a_m$ are real constants. Moreover, from (\ref{conclusao.xi.final}) and (\ref{conclusao.A2}), we have that
\begin{gather*}
A^{ij} - a - A^{ik} - A^{kj} = 0  \qquad {\rm and} \qquad	
A^{ik} - a - A^{ij} - A^{jk} = 0.	
\end{gather*}
By taking the sum of these equations and using~(\ref{a.simetrico}), we obtain $A^{kj}=-a$. Therefore, it follows from~(\ref{conclusao.fi}) and~(\ref{conclusao.B}), that
\begin{gather}
	\phi^{st} = -a h_{st}, \qquad \forall \, s \neq t.
	\label{conclusao.phi.final}
\end{gather}
Moreover, from (\ref{conclusao.xi.final}) and (\ref{eta.l}), we get
\begin{gather}
\eta^i_{,l_i} l_j = \eta^j. \label{etai.etaj}
\end{gather}
 Since the function $\eta^m$ depends only on $x_m$ and $l_m$, we conclude that,
$\eta^i_{,l_il_i} = 0$, i.e.,
\begin{gather}
	\eta^i = N^i(x_i)l_i + M^i(x_i).
	\label{forma.de.eta}
\end{gather}
Hence, it follows from (\ref{etai.etaj}) and (\ref{forma.de.eta}) that $\eta^i_{,l_i} = \eta^j_{,l_j} = N(x_i)$. Therefore, $N'(x_i) = \eta^j_{,l_jx_i} = 0$, which implies that,
\begin{gather}
	\eta^i = c l_i + M^i(x_i).
	\label{conclusao.eta.quase}
\end{gather}

Finally, we apply the prolongation ${\rm pr}^{(1)}V$ to equation (\ref{C}) to obtain
\begin{gather*}
	\ei \eta^i_{(i)} + \ej \phi^{ji}l_j + \ej h_{ji} \eta^j + \ek \phi^{ki} l_k + \ek h_{ki} \eta^k = 0,
\end{gather*}
which implies that
\begin{gather*}
\ei \eta^i_{,x_i} + \ei \eta^i_{,l_i} l_{i,x_i} - \ei \xi^i_{,x_i}l_{i,x_i} + \ej \phi^{ji}l_j + \ej h_{ji} \eta^j + \ek \phi^{ki} l_k + \ek h_{ki} \eta^k = 0.
\end{gather*}
When we substitute (\ref{subst.lii}) for  $l_{i,x_i}$ and we consider equations (\ref{conclusao.xi.final}), (\ref{conclusao.phi.final}) and  (\ref{conclusao.eta.quase}), we obtain
\begin{gather*}
\ei M^i_{,x_i} + \ej h_{ji} M^i + \ek h_{ki} M^k = 0. 
\end{gather*}
The analysis of the coef\/f\/icients of $h_{ji}$ and $h_{ki}$ enables us to conclude that $M^i = M^k = 0$, consequently, $\eta^m = c l_m$, $\forall \, 1 \leq m \leq 3$. This concludes the proof of Theorem~\ref{teo.generator}.
\end{proof}

\subsection*{Acknowledgements}

The authors were partially supported by CAPES/PROCAD and CNPq.

\LastPageEnding


\begin{thebibliography}{99}
\footnotesize\itemsep=0pt

\bibitem{barbosa}
Barbosa J.L., Ferreira W., Tenenblat K., Submanifolds of constant sectional
  curvature in pseudo-{R}iemannian manifolds, \href{http://dx.doi.org/10.1007/BF00129898}{\textit{Ann. Global Anal. Geom.}}
  \textbf{14} (1996), 381--401.

\bibitem{Cartan}
Cartan E., La d\'eformation des hypersurfaces dans l'espace conforme r\'eel \`a
  {$n \ge 5$} dimensions, \textit{Bull. Soc. Math. France} \textbf{45} (1917),
  57--121.

\bibitem{CorMilMart}
Corro A.V., Mart{\'{\i}}nez A., Mil{\'a}n F., Complete f\/lat surfaces with two
  isolated singularities in hyperbolic 3-space, \href{http://dx.doi.org/10.1016/j.jmaa.2009.12.020}{\textit{J.~Math. Anal. Appl.}}
  \textbf{366} (2010), 582--592, \href{http://arxiv.org/abs/0905.2371}{arXiv:0905.2371}.

\bibitem{ferreira}
Ferreira W., Solu\c{c}\~oes invariantes pelos grupos de simetria de Lie das
  Equa\c{c}\~oes Generalizadas Intr\'{\i}nsecas de Laplace e de sinh-Gordon
  el\'{\i}ptica e propriedades geom\'etricas das subvariedades associadas,
  Ph.D.~thesis, Universidade de Bras\'{\i}lia, 1994.

\bibitem{galvez}
G{\'a}lvez J.A., Mart{\'{\i}}nez A., Mil{\'a}n F., Flat surfaces in the
  hyperbolic {$3$}-space, \href{http://dx.doi.org/10.1007/s002080050337}{\textit{Math. Ann.}} \textbf{316} (2000), 419--435.

\bibitem{Guichard}
Guichard C., Sur les syst\`emes triplement ind\'etermin\'es et sur les
  syst\`emes triplement orthogonaux, Gauthier-Villars, Paris, 1905.

\bibitem{Jeromin2}
Hertrich-Jeromin U., Introduction to {M}\"obius dif\/ferential geometry,
  \href{http://dx.doi.org/10.1017/CBO9780511546693}{\textit{London Mathematical Society Lecture Note Series}}, Vol.~300, Cambridge
  University Press, Cambridge, 2003.

\bibitem{Jeromin1}
Hertrich-Jeromin U., On conformally f\/lat hypersurfaces and {G}uichard's nets,
  \textit{Beitr\"age Algebra Geom.} \textbf{35} (1994), 315--331.


\bibitem{suyamajer2}
Hertrich-Jeromin U., Suyama Y., Conformally f\/lat hypersurfaces with
  {B}ianchi-type {G}uichard nets, \textit{Osaka~J. Math.}, {t}o appear.

\bibitem{suyamajer}
Hertrich-Jeromin U., Suyama Y., Conformally f\/lat hypersurfaces with cyclic
  {G}uichard net, \href{http://dx.doi.org/10.1142/S0129167X07004138}{\textit{Internat.~J. Math.}} \textbf{18} (2007), 301--329.

\bibitem{flatfront1}
Kokubu M., Rossman W., Saji K., Umehara M., Yamada K., Singularities of f\/lat
  fronts in hyperbolic space, \href{http://dx.doi.org/10.2140/pjm.2005.221.303}{\textit{Pacific~J. Math.}} \textbf{221} (2005),
  303--351, \href{http://arxiv.org/abs/math.DG/0401110}{math.DG/0401110}.

\bibitem{flatfront2}
Kokubu M., Umehara M., Yamada K., Flat fronts in hyperbolic 3-space,
  \href{http://dx.doi.org/10.2140/pjm.2004.216.149}{\textit{Pacific~J. Math.}} \textbf{216} (2004), 149--175,
  \href{http://arxiv.org/abs/math.DG/0301224}{math.DG/0301224}.

\bibitem{Lafontaine}
Lafontaine J., Conformal geometry from the {R}iemannian viewpoint, in Conformal
  Geometry ({B}onn, 1985/1986), \textit{Aspects Math.}, Vol.~E12, Vieweg,
  Braunschweig, 1988, 65--92.

\bibitem{reflame}
Lam\'e G., Le\c{c}ons sur les coordonn\'es curvilignes et leurs diverses
  applications, Mallet-Bachelier, Paris, 1859.

\bibitem{lie}
Lie S., Theorie der Transformationsgruppen, B.G.~Teubner, Leipzig, 1888, 1890,
  1893.

\bibitem{marsanten}
Martinez A., dos Santos J.P., Tenenblat K., Helicoidal f\/lat surfaces in the
  hyperbolic 3-space, \textit{Pacific~J. Math.}, {t}o appear.

\bibitem{Olver}
Olver P.J., Applications of {L}ie groups to dif\/ferential equations,
  \href{http://dx.doi.org/10.1007/978-1-4684-0274-2}{\textit{Graduate Texts in Mathematics}}, Vol.~107, Springer-Verlag, New York,
  1986.

\bibitem{olver2}
Olver P.J., Symmetry groups and group invariant solutions of partial
  dif\/ferential equations, \textit{J.~Differential Geom.} \textbf{14} (1979),
  497--542.


\bibitem{rabelo}
Rabelo M.L., Tenenblat K., Submanifolds of constant nonpositive curvature,
  \textit{Mat. Contemp.} \textbf{1} (1991), 71--81.

\bibitem{Suyama1}
Suyama Y., Conformally f\/lat hypersurfaces in {E}uclidean 4-space,
  \textit{Nagoya Math.~J.} \textbf{158} (2000), 1--42.

\bibitem{Suyama2}
Suyama Y., Conformally f\/lat hypersurfaces in {E}uclidean 4-space.~{II},
  \textit{Osaka~J. Math.} \textbf{42} (2005), 573--598.

\bibitem{Suyama3}
Suyama Y., Conformally f\/lat hypersurfaces in {E}uclidean 4-space and a class of
  {R}iemannian 3-manifolds, \textit{S\=urikaisekikenky\=usho K\=oky\=uroku}
  (2001), no.~1236, 60--89.


\bibitem{keti}
Tenenblat K., Transformations of manifolds and applications to dif\/ferential
  equations, \textit{Pitman Monographs and Surveys in Pure and Applied
  Mathematics}, Vol.~93, Longman, Harlow, 1998.

\bibitem{keti3}
Tenenblat K., Winternitz P., On the symmetry groups of the intrinsic
  generalized wave and sine-{G}ordon equations, \href{http://dx.doi.org/10.1063/1.530042}{\textit{J.~Math. Phys.}}
  \textbf{34} (1993), 3527--3542.

\end{thebibliography}
\end{document}